\documentclass{article}

\usepackage{amsmath,amsfonts,amssymb,amsthm,pstricks,pst-plot,pst-node,eucal}
\usepackage[dvips]{graphics}

\newcommand{\Z}{\mathbb{Z}}
\newcommand{\R}{\mathbb{R}}
\newcommand{\N}{\mathbb{N}}
\newcommand{\oh}{\tfrac12}

\DeclareMathOperator{\Prob}{Prob}

\DeclareMathOperator{\Ai}{Ai}

\DeclareMathOperator{\Li}{Li}

\newcommand{\cl}{\{\lambda(t)\}}

\newcommand{\lang}{\left\langle}
\newcommand{\rang}{\right\rangle}
\newcommand{\gl}{\mathfrak{gl}}

\newpsobject{showgrid}{psgrid}{subgriddiv=1,griddots=5,gridlabels=0pt}

\newtheorem{theorem}{Theorem}

\newcommand{\LV}{\Lambda^{\frac\infty2}V}
\newcommand{\ul}{\underline}
\newcommand{\al}{\alpha}
\newcommand{\e}{\varepsilon}
\newcommand{\zp}{z\frac{\partial}{\partial z}}

\theoremstyle{definition}
\newtheorem{remark}{Remark}
\newtheorem{definition}{Definition}

\begin{document}

\title
{Random skew plane partitions and the Pearcey process}
\author{Andrei Okounkov \thanks{Department of Mathematics, Princeton
    University, Princetion NJ 08544-1000. E-mail: okounkov@math.princetion.edu}
and Nicolai Reshetikhin\thanks{Department of Mathematics, University of California at
Berkeley, 
Berkeley, CA 94720-3840. E-mail:  
reshetik@math.berkeley.edu }}
\date{}
\maketitle

\begin{abstract}
We study random skew 3D partitions weighted by $q^{\textup{vol}}$ and, 
specifically, the $q\to 1$ asymptotics of local correlations near various
points of the limit shape. We obtain sine-kernel asymptotics for 
correlations in the bulk of the disordered region, Airy kernel 
asymptotics near a general point of the frozen boundary, and 
a Pearcey kernel asymptotics near a cusp of the frozen 
boundary.  
\end{abstract}

\section{Introduction}

\begin{figure}[!htbp]
  \centering
  \scalebox{0.33}{\includegraphics{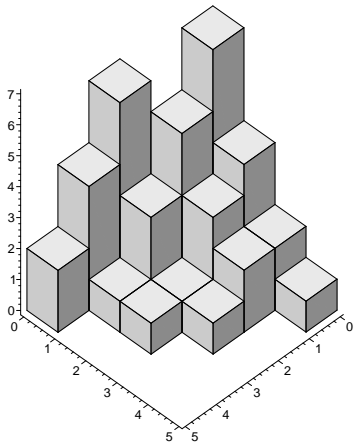}}
  \caption{A skew 3D partition}
  \label{fig_skew}
\end{figure}

A plane partition $\pi=(\pi_{ij})$ is an array of nonnegative numbers
indexed by $(i,j)\in \N^2$ that is monotone, that is, 
$$
\pi_{ij} \ge \pi_{i+r,j+s}\,, \quad r,s\ge 0 
$$
and  finite in the sense that $\pi_{ij}=0$ when $i+j\gg 0$. 
Plane partitions have a obvious generalization which we 
call \emph{skew plane partitions}. A skew plane partition is again 
a monotone array $(\pi_{ij})$ which is now indexed by points 
$(i,j)$ of a skew shape $\lambda/\mu$, where $\mu\subset\lambda$ is a pair
of ordinary partitions. We call $\mu$ and $\lambda$ the 
\emph{inner} and \emph{outer} shape of $\pi$, respectively. 
In fact, in this paper we will only consider the case when 
the outer shape $\lambda$ is a $a\times b$ rectangle. 
Here is an example with $\mu=(1,1)$ and $a=b=5$.  
\begin{equation}
  \label{pi}
 \pi =
\begin {array}{ccccc}
&7&4&2&1\\
&5&3&2&0\\
6&3&1&1&0\\
4&1&1&0&0\\
2&0&0&0&0
\end {array}
\quad \,. 
 \end{equation}
Placing $\pi_{ij}$ cubes over the $(i,j)$ square in $\lambda/\mu$
gives a three-dimensional object which we will call a 
\emph{skew 3D partition} and denote by the same letter $\pi$. 
Its volume is $|\pi|=\sum \pi_{ij}$. 
For $\pi$ as in \eqref{pi}, it shown in Figure \ref{fig_skew}. 

Given a parameter $0<q<1$, define a probability measure on 
the set of all skew plane partitions with given inner and outer 
shapes by setting
\begin{equation}\label{distr_}
\Prob(\pi)\propto q^{|\pi|} \,.
\end{equation}
The corresponding random skew 3D partition model has 
a natural random growth interpretation, the parameter
$q$ being the fugacity. Also, a simple bijection, which 
should be clear from Figure \ref{fig_skew} and is 
recalled below, relates this model to a random tiling 
problem. 

\begin{figure}[!htbp]
  \begin{center}
    {\scalebox{0.4}{\includegraphics*[50,110][520,650]{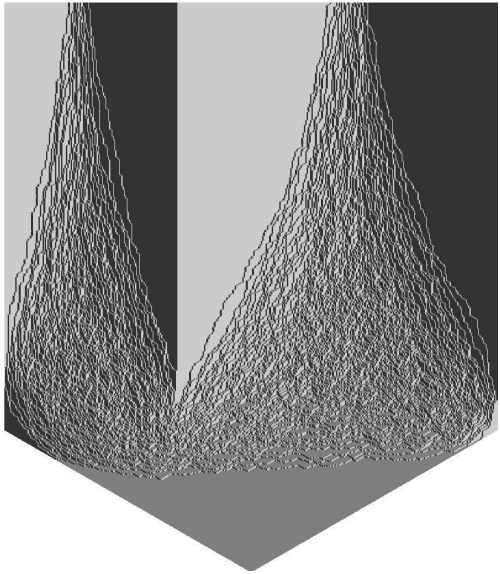}}}
   \caption{A large random skew 3D partition}
    \label{skew-random}
  \end{center}
\end{figure}

We are interested in the thermodynamic limit in which 
$q\to 1$ and both inner and outer shapes are rescaled by $1/r$
where
$$
r=-\ln q \to + 0 \,.  
$$
The results of \cite{CKP} imply the following form of the law of 
large numbers: scaled by $r$ in all directions, the surface of 
our random skew 3D partition converges to a nonrandom surface ---
the limit shape. This limit shape will be easy to see in the exact
formulas discussed below. A simulation showing the formation 
of the limit shape is presented in Figure \ref{skew-random}. 

An important qualitative feature of limit shape is 
the presence of both ordered and disordered regions, 
separated by the \emph{frozen boundary}.  Furthermore, 
the frozen boundary has various special points, namely, 
it has \emph{cusps} (there is one forming in Figure \ref{skew-random}, it 
can be seen more clearly in Figures \ref{shape1} and 
\ref{shape2}) and also 
\emph{turning points} where the limit shape is 
not smooth. In Figure \ref{bc1}--\ref{shape2}, the turning points 
are the points of tangency to 
any of the lines in the same figure. 

One expects that the microscopic properties of the random surface, 
in particular, the correlation functions of local operators, are 
\emph{universal} in the sense that they are determined by 
the macroscopic behavior of the limit shape at that point.  
More specifically, one expects that:  
\begin{enumerate}
\item in the bulk of the disordered region, the correlation are 
given by the incomplete beta kernel \cite{K,OR} with the 
parameters determined by slope of the limit shape (a special 
case of this is the discrete sine kernel);
\item at a general point of the frozen boundary, suitably scaled, 
the correlation are given by the extended Airy kernel \cite{PS,J1,TW1};
\item at a cusp of the frozen boundary, correlation, suitably
scaled, are given by the extended Pearcey kernel, discussed 
below and in  \cite{TW2}.  
\end{enumerate}
In this paper we prove all these statements for the model 
at hand. The required techniques were developed in our paper \cite{OR}, 
of which this one is a continuation. Namely, as will be 
reviewed below, our random skew
plane partition model is a special case of \emph{Schur process}. 
This yields an exact contour integral formulas for correlation 
functions. The asymptotics is then extracted by a direct
albeit laborious saddle point analysis. The striking resemblance 
of the above list to classification of singularities is not 
accidental for, as we will see, these three situation correspond 
precisely to the saddle point being a 
 simple, double, or triple critical point.  

We will also see that the 
frozen boundary is essentially an algebraic curve and 
that it has precisely one cusp per each exterior corner 
of the inner shape $\mu$.  
 
We expect that near a turning point the correlations 
behave like eigenvalues of a $k\times k$ corner 
of a GUE random $N\times N$ matrix, where $N\gg 0$ and 
$k$ plays the role of time. We hope to return to this 
question in a future paper.

The results presented here were obtained in 2002-03 and were reported 
by us at several conferences. The period between then and now saw
many further developments in the field. Most notably, the 
Pearcey process, which we found describes the behavior near a cusp 
of the frozen boundary arose in the random matrix context in the 
work of Tracy and Widom \cite{TW2}. Pearcey asymptotics for equal time 
correlations of eigenvalues were obtained earlier by Brezin and Hikami 
\cite{BH1,BH2}
and also by Aptekarev, Bleher, and Kuijlaars \cite{ABK}. We enjoyed and 
benefited from the correspondence with C.~Tracy on 
subject. 

In \cite{FS},
Ferrari and Spohn derived from the exact formulas of \cite{OR} the 
Airy process asymptotics in the case of unrestricted 3D partitions
(the $\mu=\emptyset$, $a=b=\infty$ case in our notation). In a 
related but technically more involved 
context, the Airy process asymptotics was found by 
K.~Johansson in \cite{J2}. 

In \cite{ORV}, the partition function of the 
random surface model studied here was related to 
the \emph{topological vertex} of \cite{AKMV} and, 
thus, to the Gromov-Witten theory 
of toric Calabi-Yau threefolds. They were many 
subsequent developments, some of which
are reviewed in \cite{O}. The papers \cite{DST,SV} may 
the closest to the material presented here. 
 Also, much more general results on algebraicity 
of the frozen boundary are now available \cite{KO}. We would like 
to thank R.~Kenyon and C.~Vafa for numerous discussions.  

N.~R.\ is grateful to Laboratoire de Physique
Theorique at Saclay for the hospitality, where
part of this work was done and to J.-B.\ Zuber and Ph.\ Di Francesco
for interesting discussions. His work was supported by the NSF grant
DMS-0070931 and by the Humboldt Foundation. The work of A.~O.\ was 
partially supported by the Packard Foundation.

\section{Preliminaries}\label{comb}
\subsection{Skew 3D partition as a sequence of its slices}

We associate to $\pi$ the sequence $\cl$ of its diagonal slices,
that is, the sequence of partitions
\begin{equation}
  \label{pila}
 \lambda(t)=(\pi_{i,t+i})\,, \quad i\ge\max(0,-t), \ t\in {\mathbb Z} \,.
\end{equation}
Throughout this paper, we assume that the outer shape of our 
skew partition is an $a\times b$ box and, in particular, we 
will use the letter $\lambda$ to denote diagonal slices, not 
the outer shape.

Notation $\lambda \succ \nu$ as usual
means that $\lambda$ and $\nu$ interlace,
that is,
$$
\lambda_1 \ge \mu_1 \ge \lambda_2 \ge \mu_2 \ge \lambda_3 \ge \dots \,.
$$
It is easy to see that the sequence $\cl$ corresponds
to a skew plane partition if and only if it satisfies the
following conditions:
\begin{itemize}
\item if the slice $\lambda(t_0)$ is passing through an inner
corner of the skew plane partition then 
\begin{equation}
  \label{laint}
 \dots \prec \lambda(t_0-2) \prec \lambda(t_0-1) \prec \lambda(t_0)
\succ \lambda(t_0+1) \succ \lambda(t_0+2) \succ \dots\,; 
\end{equation}
\item if the slice $\lambda(t_0)$ is passing through an
outer corner of the skew plane partition then 
\begin{equation}
  \label{inner}
 \dots \succ \lambda(t_0-2) \succ \lambda(t_0-1) \succ \lambda(t_0)
\prec \lambda(t_0+1) \prec \lambda(t_0+2) \prec \dots\,. 
\end{equation}
\end{itemize}

\noindent
For example, the configuration $\cl$ corresponding to
the partition \eqref{pi} is
\begin{equation*}
(2)\prec (4)\prec (6,1)\succ (3,1) \prec (5,1) \prec 
(7,3,1) \succ (4,2) \succ (2) \succ(1) \,.
\end{equation*}
We will denote the sequence of inner and outer corners of the inner shape 
by $\{v_i\}_{1\leq i
\leq N}$ and $\{u_i\}_{1\leq i\leq
N-1}$, respectively. We assume that they are numbered so that 
$$
v_1<u_1<v_2<u_2<\dots u_{N-1}<v_N \,. 
$$
We also assume that the point $t=0$ is chosen so that 
\begin{equation}\label{norm-uv}
\sum_{1\leq
i\leq  N}v_i=\sum_{1\leq i\leq  N-1}u_i \,.
\end{equation}

\subsection{Connection to tilings}

\begin{figure}[htbp]\psset{unit=0.5cm}
  \begin{center}
    \begin{pspicture}(-5,-5)(5,7)
\scriptsize \showgrid 
\psset{dimen=middle}
\psaxes[axesstyle=frame,Ox=-5,Oy=-5,Dx=2,Dy=2,ticks=none](-5,-5)(5,7)
\psdiamond[fillstyle=solid,fillcolor=lightgray](-4,-.5)(1,0.5)
\psdiamond[fillstyle=solid,fillcolor=lightgray](-3,-3)(1,0.5)
\psdiamond[fillstyle=solid,fillcolor=lightgray](-3,2)(1,0.5)
\psdiamond[fillstyle=solid,fillcolor=lightgray](-2,-3.5)(1,0.5)
\psdiamond[fillstyle=solid,fillcolor=lightgray](-2,-1.5)(1,0.5)
\psdiamond[fillstyle=solid,fillcolor=lightgray](-2,4.5)(1,0.5)
\psdiamond[fillstyle=solid,fillcolor=lightgray](-1,-4)(1,0.5)
\psdiamond[fillstyle=solid,fillcolor=lightgray](-1,-2)(1,0.5)
\psdiamond[fillstyle=solid,fillcolor=lightgray](-1,1)(1,0.5)
\psdiamond[fillstyle=solid,fillcolor=lightgray](0,-4.5)(1,0.5)
\psdiamond[fillstyle=solid,fillcolor=lightgray](0,-3.5)(1,0.5)
\psdiamond[fillstyle=solid,fillcolor=lightgray](0,-1.5)(1,0.5)
\psdiamond[fillstyle=solid,fillcolor=lightgray](0,3.5)(1,0.5)
\psdiamond[fillstyle=solid,fillcolor=lightgray](1,-4)(1,0.5)
\psdiamond[fillstyle=solid,fillcolor=lightgray](1,-2)(1,0.5)
\psdiamond[fillstyle=solid,fillcolor=lightgray](1,1)(1,0.5)
\psdiamond[fillstyle=solid,fillcolor=lightgray](1,6)(1,0.5)
\psdiamond[fillstyle=solid,fillcolor=lightgray](2,-3.5)(1,0.5)
\psdiamond[fillstyle=solid,fillcolor=lightgray](2,-.5)(1,0.5)
\psdiamond[fillstyle=solid,fillcolor=lightgray](2,2.5)(1,0.5)
\psdiamond[fillstyle=solid,fillcolor=lightgray](3,-3)(1,0.5)
\psdiamond[fillstyle=solid,fillcolor=lightgray](3,0)(1,0.5)
\psdiamond[fillstyle=solid,fillcolor=lightgray](4,-1.5)(1,0.5)
\end{pspicture}
    \caption{Horizontal tiles of the tiling corresponding to the partition 
    in Figure \ref{fig_skew} in $(t,h)$-coordinates}
    \label{f3}
  \end{center}
\end{figure}
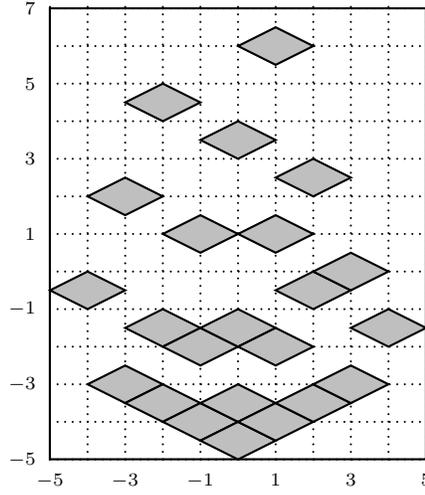

There is a well-known mapping of 3D diagrams to tilings of the
plane by rhombi. Namely, the tiles
are the images of faces of the 3D diagram under the projection
\begin{equation}
(x,y,z)\mapsto (t,h)=(y-x,z-(x+y)/2)\,.\label{xyz}
\end{equation}
This mapping is a bijection between 3D diagrams and tilings with
appropriate boundary conditions. The horizontal tiles of the 
tiling corresponding to the
diagram in Figure \ref{fig_skew} are shown in Figure \ref{f3}.

It is clear that the positions of
horizontal tiles uniquely determine both the tiling and the
partition $\pi$. The set
\begin{equation}
  \label{def_sigm}
    \sigma(\pi)= \{(i-j,\pi_{ij}-(i+j-1)/2)\}\subset \Z \times \oh \, \Z
\end{equation}
is precisely the set of the centers of the horizontal tiles.
Notice that if $(h,t)$ is a center of a tile $h+t/2+1/2$ is always
an integer.

Define
$$
B(t)= \frac12\sum_{i=1}^N \left| t - v_i \right| -
\frac12\sum_{i=1}^{N-1} \left| t - u_i \right| \,.
$$

\begin{figure}[htbp]
  \begin{center}
    \scalebox{0.3}{\includegraphics{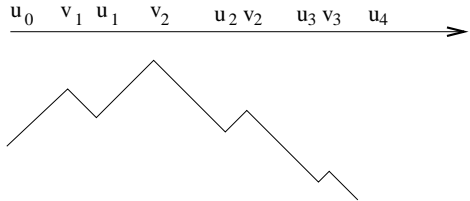}}
    \caption{Coordinates of corners}
    \label{uv.eps}
  \end{center}
\end{figure}

The image of the inner boundary of our skew plane partitions  in
the $(h,t)$-plane is the curve
\[
h=-B(t)\,,
\]
see an example of this curve in Figure \ref{uv.eps}. 
In particular, the highest layer of horizontal
rhombi for an empty plane partitions is the set of points with
coordinates $h=-B(t)-1/2$.

\subsection{Partition function and correlation functions}

Generalizing \eqref{distr_}, introduce a probability measure 
on skew plane partitions by 
\begin{equation}\label{epschur}
\Prob(\{\lambda(t)\}) \,\propto\,
\prod_{t\in \Z} q_{t}^{|\lambda(t)|}\,,
\end{equation}
where $0\leq q_t <1$ are parameters. 

We assume that $q_t=0, \ t<u_0=-a$ or $t>u_N=b$ and so 
the plane partition is confined to a  $a\times b$ outer box. 
The homogeneous case when all nonzero $q_t$ are equal corresponds
to \eqref{distr_}. 

For fixed inner shape $\mu$, the partition function is 
defined by 
\[
Z(\{q_t\},\mu)=\sum_{\cl}
\prod_{t\in \Z} q_{t}^{|\lambda(t)|}=\sum_\pi\prod_{t\in \Z}q_t^{|\pi_t|}\, .
\]
where $|\pi_t|=\sum_i \pi_{i,t+i}$.

The correspondence 
$\pi\mapsto \sigma(\pi)$ defined in \eqref{def_sigm} makes a 
random skew partition a random subset of $\Z\times
(\Z+\oh)$, that is, a random point field on a lattice. 
This motivates the following 

\begin{definition}
 Given a subset $U\subset \Z\times (\Z+\oh)$, define
the corresponding correlation function by
\begin{equation}\label{cor}
\rho(U)=\Prob\left(U\subset \sigma(\pi)\right)
=\frac{1}{Z}\sum_{\pi, U\subset \sigma(\pi)} \prod_{t\in \Z}q_t^{|\pi_t|}\,.
\end{equation}
These correlation functions depend
on parameters $q_t$ and on the fixed inner shape $\mu$ of skew
plane partitions.
\end{definition}

Consider the following ``local''  functions  on skew
plane partitions :
\[
\rho_{h,t}(\pi)=\left\{\begin{array}{ll}1& \mbox{if} \
(h,t)\in \sigma(\pi),\\
0 & \mbox{otherwise} .
\end{array}\right.
\]

If $U=\{(h_1,t_1),\dots, (h_n,t_n)\}$ with $t_1\geq \dots \geq t_n$
and $(h_i,t_i)\neq (h_j,t_j)$, the
correlation function (\ref{cor}) can be written as:
\begin{equation}\label{cor1}
\rho(U)=\lang\rho_{h_1,t_1}\dots \rho_{h_n,t_n}\rang=\frac{1}{Z}\sum_{\pi}
\rho_{h_1,t_1}(\pi)\dots\rho_{h_n,t_n}(\pi)\prod_{t\in\Z}q_t^{|\pi_t|} \,. 
\end{equation}

\section{Schur processes}\label{schur}

\subsection{General Schur processes}
Schur process, introduced in \cite{OR} is a 
probability measure on
sequences of partitions. 

Parameters of the Schur process are
sequences of pairs of functions $\{\phi_t^\pm(z)\}_{t\in \Z}$
such that $\phi^+(z)$ is analytic at $z=0$ and $\phi^-(z)$ is
analytic at $z=\infty$.
For such  pair of functions $\phi(z)^\pm$  consider
skew Schur functions
\[
s_{\lambda/\mu}[\phi^+]=\det(\phi^+_{\lambda_i-\mu_j-i+j}) \ ,
\]
and
\[
s_{\lambda/\mu}[\phi^-]=\det(\phi^-_{-\lambda_i+\mu_j+i-j}) \ .
\]
Some basic notions about Schur functions are recalled in the 
Appendix \ref{app_schur}. 

Define the transition weight by the formula
\[
S_\phi(\lambda,\mu)=\sum_\nu s_{\lambda/\nu}[\phi^+]s_{\mu/\nu}[\phi^-]
\]

\begin{definition}\label{defS}
The
probabilities of the Schur process are given by
$$
\Prob(\{\lambda(t)\}) = \frac1Z \, \prod_{m\in \Z+1/2}
S_{\phi_m}\left(\lambda(m-\oh),\lambda(m+\oh)\right)\,,
$$
where the transition weight $S_\phi$ is defined above, by $\phi_t$
we denoted the pair of functions $\phi_t^\pm$ and $Z$ is the
normalizing factor (partition function)
$$
Z=\sum_{\cl} \prod_{m\in \Z+1/2}
S_{\phi[m]}\left(\lambda(m-\oh),\lambda(m+\oh)\right) \,.
$$
\end{definition}

If instead of infinite sequences $\cl$ we have  finite  sequences of
length $N$ 
we will say that the Schur
process is of length $N$.

\subsection{Polynomial Schur processes and height distributions on
skew plane partitions}
We will say that the Schur process is {\it polynomial} if functions
$\phi^\pm_t(z)$ are polynomials in $z^\pm$.
Let us show that the measure (\ref{epschur}) is closely
related to  a polynomial Schur process.

We will parametrize the inner shape as before by
assuming that $v_i, \ 1\leq i\leq N , \ v_i\in \Z$
are positions of inner corners of the inner shape of the plane
partition (see Fig. \ref{uv.eps}) and $u_i\in (v_i, v_{i+1}), \
1\leq i\leq N-1$
are positions of outer corners.

\begin{theorem}The restriction of the measure (\ref{epschur})
to random variables supported on subsequences $\{\lambda(m)\}
, \ m<v_1$, \ $\{\lambda(v_i)\}, \ 1\leq i
\leq N$, \ $\{\lambda(m)\}, \ m>v_N$ coincides with the polynomial Schur
process with parameters
\begin{eqnarray}\label{x-pm}
\phi_m^-(z)&=&(z-x^-_m)\,,\quad \phi_m^+(z)=1, \ \\
\phi_i^+(z)&=&\prod_{v_i<m<u_i, \ m\in \Z+\oh}(z-x^+_m) \ , \\
\phi_i^-(z)&=&\prod_{u_i<m<v_{i+1}, \ m\in \Z+\oh}(z-x^-_m) \ ,\\
\phi_m^+(z)&=&(z-x^+_m)\,, \quad \phi_m^-(z)=1 , \
\end{eqnarray}
where parameters $x^\pm_t$ and $q_t$ are related as follows:
\begin{eqnarray}\label{q-x}
\frac{x^+_{m+1}}{x^+_m}&=&q_{m+\oh}, \  v_i<m<u_i-1, \mbox{or} \ m>v_N\\
x^+_{u_i-\oh}x^-_{u_i+\oh}&=&q_{u_i}^{-1}, \\
x^-_{v_i-\oh}x^+_{v_i+\oh}&=&q_{v_i}, \\
\frac{x^-_m}{x^-_{m+1}}&=&q_{m+\oh},  u_i<m<v_{i+1}-1, \mbox{or} \
m<v_1
\end{eqnarray}
\end{theorem}
\begin{proof}Let us restrict the process (\ref{epschur}) to the
subsequence $\{\lambda_{v_i}\}$.
It is easy to see that transition probability from $v_{i+1}$ to $v_i$ in
such subprocess are
\begin{eqnarray}
S^i(\lambda(v_i), \lambda(v_{i+1})&=\sum s_{\lambda(v_i)/\lambda(u_i)}
(x^+_{v_i+1/2},x^+_{v_i+3/2},\dots,x^+_{u_i-1/2}) \\
&s_{\lambda(v_{i+1})/\lambda(u_i)}
(x^-_{v_{i+1}-1/2},x^-_{v_{i+1}-3/2},\dots ,x^-_{u_i+1/2})
\end{eqnarray}
where $x^\pm_m$ are related to $q_t$ as in (\ref{q-x}). Thus this process
(\ref{epschur}) is a polynomial Schur process with $\phi^\pm_i$ given
by (\ref{x-pm}).
Conversely, it is clear, that due to the identity
(\ref{recurschur}) any polynomial Schur
process can be extended to a probability measure (\ref{epschur})
on sequences of interlacing partitions with parameters
$\{q_t\}$ defined as in (\ref{x-pm})(\ref{q-x}).
\end{proof}

\section{Fermionic representation for correlation functions}\label{ferm}

\subsection{}

For  $m\in \Z + \oh$ define $\e(m)=+$ if $v_i<m<u_i$ and $1\leq
i\leq N$, and $\e(m)=-$ if $u_i<m<v_{i+1}$ and $0\leq i\leq N-1$.
This is shown on Fig.\ref{uv.eps}. Define $ D^+=\{m|\e(m)=+\}$ and
$D^-=\{m|\e(m)=-\}$.

Let $x^{\pm}_m$ be positive numbers related to $q_t$ as in
(\ref{q-x}). 
Notice that for given $q_t$ the numbers $x^\pm_m$ are
defined up to a transformation $x^\pm_m\to x^\pm_m a^{\pm 1}$.

\begin{theorem}\begin{enumerate}
\item The partition function for the height distribution on skew
plane partitions can be represented as the matrix element of the
product of vertex operator described in the Appendix \ref{vo} as
follows
\begin{multline}
Z=\Bigg(\prod_{u_N>m>v_N}\Gamma_-(x^+_m)\dots
\prod_{u_i<m<v_{i+1}}\Gamma_+(x^-_m)
\prod_{v_i<m<u_{i}}\Gamma_-(x^+_m)
\\
\prod_{v_1>m>u_0}\Gamma_+(x^-_m)v^{(m)}_0,
v^{(m)}_0\Bigg)= \Bigg(\prod_{m\in \Z+\oh,
u_0<m<u_N}\Gamma_{-\e(m)}(x^{\e(m)}_m)v^{(0)}_0, v^{(0)}_0\Bigg)
\end{multline}
and
\[
Z=\prod_{m_1<m_2, m_1\in D^-, m_2\in D^+}(1-x^-_{m_1}x^+_{m_2})^{-1},\  m_i\in
\Z+\oh
\]
\item Assume $t_1 >\dots > t_n$, then
\begin{multline}\label{corr}
\lang\rho_{h_1,t_1}\dots \rho_{h_n,t_n}\rang=\frac{1}{Z}\Bigg(\prod_{m\in
\Z+\oh, u_N>m>t_1}
\Gamma_{-\e(m)}(x^{\e(m)}_m)\psi_{j_1}\psi^*_{j_1}\dots\\
\prod_{t_i<m<t_{i-1}}\Gamma_{-\e(m)}(x^{\e(m)}_m)\psi_{j_i}\psi^*_{j_i}
\prod_{t_{i+1}<m<t_i}\Gamma_{-\e(m)}(x^{\e(m)}_m)\\
\dots \psi_{j_n}\psi^*_{j_n} \prod_{t_n>m>u_0}
\Gamma_{-\e(m)}(x^{\e(m)}_m)v^{(0)}_0,v^{(0)}_0\Bigg)
\end{multline}
Here and below $j_i=h_i-B(t_i)+1/2$.
\item Correlation functions
(\ref{cor1}) are determinants:
\begin{equation}\label{main-corr}
\lang\rho_{h_1,t_1}\dots \rho_{h_n,t_n}\rang=
\det(K(t_i,h_i),(t_k,h_k)))_{1\leq i,k\leq n}
\end{equation}
where
\begin{multline}\label{main-corr2}
K((t_1,h_1),(t_2,h_2))=\\ \frac{1}{(2\pi i)^2}
\int_{|z|<R(t_1)}\int_{R^*(t_2)<|w|}
\frac{\Phi_-(z,t_1)\Phi_+(w,t_2)}{\Phi_+(z,t_1)\Phi_-(w,t_2)}\\
\frac{\sqrt{zw}}{z-w}z^{-h_1+B(t_1)-1/2}w^{h_2-B(t_2)+1/2}\frac{dzdw}{zw}
\end{multline}
Here $|w|<|z|$ for $t_1\geq t_2$, $|w|>|z|$ for $t_1<t_2$,
$R(t)=\min_{m>t}(|x^+_m|^{-1})$ and $R^*(t)=\max_{m<t}(|x^-_m|)$. Functions
$\Phi_\pm(z,t)$ are:
\begin{eqnarray}
\Phi_+(z,t)&=&\prod_{m>t, m\in D^+, m\in \Z+\oh}(1-zx_m^+)\\
\Phi_-(z,t)&=&\prod_{m<t, m\in D^-, m\in \Z+\oh}(1-z^{-1}x_m^-)
\end{eqnarray}
\end{enumerate}
\end{theorem}
\begin{proof}

The fact that the partition function and correlation functions for
the height distribution of plane partitions the  matrix element of
the product of vertex operators as above follows form the formula
(\ref{me}) for matrix elements of products of vertex operators
$\Gamma_\pm(x)$.

Using the commutation relations (\ref{gg}), and the fact that
$\Gamma_-(x)v^(m)_0=0$ we obtain the product formula for the
partition function.

The operators $\psi_j\psi_j^*$ act on the vector $v_\lambda^{(0)}$
as follows:
\[
\psi_j\psi_j^*v_\lambda^{(0)}=v_\lambda^{(0)}
\]
if $j=\lambda_i-i+1/2$ for some $i=1,2,\dots$ and
\[
\psi_j\psi_j^*v_\lambda^{(0)}=0
\]
otherwise. Using this fact and the formula for the matrix elements
of $\Gamma_\pm(x)$ we obtain the formula (\ref{corr}) for the
correlation functions of densities.

Moving operators $\Gamma_-$ to the right and $\Gamma_+$ to the
left  we obtain the following formula for the correlation
functions
\begin{multline*}
\lang\rho_{h_1,t_1}\dots \rho_{h_n,t_n}\rang=\\
(\psi_{j_1}(t_1)\psi^*_{j_1}(t_1)\dots\psi_{j_i}(t_i)\psi^*_{j_i}(t_i)
\dots  \psi_{j_n}(t_n)\psi^*_{j_n}(t_n)v^{(0)}_0,v^{(0)}_0)
\end{multline*}
where
\begin{multline}\label{psi}
\psi_j(t)=\prod_{m>t, m\in D^+}\Gamma_-(x^+_m)\prod_{m<t, m\in
D^-} \Gamma_+(x^-_m)^{-1}\psi_j \\
 \prod_{m<t, m\in D^-} \Gamma_+(x^-_m)\prod_{m>t, m\in
D^+}\Gamma_-(x^+_m)^{-1}
\end{multline}
and
\begin{multline}\label{psi-star}
\psi^*_j(t)=\prod_{m>t, m\in D^+}\Gamma_-(x^+_m)\prod_{m<t, m\in
D^-} \Gamma_+(x^-_m)^{-1}\psi_j^* \\
 \prod_{m<t, m\in D^-} \Gamma_+(x^-_m)\prod_{m>t, m\in
D^+}\Gamma_-(x^+_m)^{-1}
\end{multline}

Here the operators on the right are given by power series.
Commuting formal power series gives the following
identities:
\[
\Gamma_+(x)^{-1}\psi_k\Gamma_+(x)=\psi_k-x\psi_{k+1}
\]
\[
\Gamma_-(x)\psi_k\Gamma_-(x)^{-1}=\sum_{n\geq 0}x^n\psi_{k-n}
\]
\[
\Gamma_+(x)^{-1}\psi^*_k\Gamma_+(x)=\sum_{n\geq 0} x^n\psi^*_{k-n}
\]
\[
\Gamma_-(x)\psi^*_k\Gamma_-(x)^{-1}=\psi^*_k-x\psi^*_{k+1}
\]



Applying these identities to the formal Fourier transform
of $\psi_j(t)$ and $\psi^*_j(t)$ we obtain:
\begin{equation}\label{ps}
\psi(z,t)=\prod_{m>t, m\in D^+}(1-zx_m^+)^{-1}\prod_{m<t,
m\in D^-}(1-z^{-1}x_m^-)
\psi(z) \ ,
\end{equation}
\begin{equation}\label{ps-star}
\psi^*(z,t)=\prod_{m>t, m\in D^+}(1-zx_m^+)\prod_{m<t, m\in D^-}
(1-z^{-1}x_m^-)^{-1} \psi^*(z)\,,
\end{equation}
where both sides are power series in $x^\pm_m$ and are
formal Laurent power series in $z$. If $v\in F$, then
$\psi(z)v\in z^{1/2}F[z^{-1}, z]]$ and
$\psi^*(z)v=z^{1/2}F[z^{-1}, z]]$.

For the inverse Fourier transform of $\psi(z,t)v$ and $\psi(z,t)$
we obtain the following
integral representations:
\begin{equation}\label{p}
\psi_k(t)v=\frac{1}{2\pi
i}\int_{|z|<R(t)}\frac{\Phi_-(z,t)}{\Phi_+(z,t)}z^{-k-1}\psi(z)vdz
\end{equation}
\begin{equation}\label{p-star}
\psi^*_k(t)v=\frac{1}{2\pi
i}\int_{R^*(t)<|w|<1}\frac{\Phi_+(w,t)}{\Phi_-(w,t)}w^{k-1}\psi^*(w)vdw
\end{equation}
Here $v\in F$, both sides are vectors in $F[[x^\pm_m]]$ and these
power series converge for sufficiently small $x$'s.

The contour of integration for $z$ is
chosen in such a way that none of the poles of $\Phi_+(z,t)$ will
be inside of it, this gives $|z|<R(t)=\min_{m>t}(|x^+_m|^{-1})$. The
contour of integration for $w$ is such that none of the poles of
$\Phi_-(w,t)$ are outside the contour. This gives
$|w|>R^*(t)=\max_{m<t}(|x^-_m|)$.


As it follows from the Wick's lemma (\ref{wl}) that correlation
functions (\ref{cor1}) are determinants of matrices of correlation
functions of two Clifford operators (\ref{psi})(\ref{psi-star}).
\[
\lang\rho_{h_1,t_1}\dots \rho_{h_n,t_n}\rang= \det(K_{a,b})_{1\leq a,b\leq
n}
\]
where
\[
K_{a,b}=K((t_a,h_a),(t_b,h_b))=
\begin{cases}
  \left(\psi_{j_a}(t_a)\, \psi^*_{j_b}(t_b) \, v_0,v_0\right)\,,
& a \le b \,,\\
-\left(\psi^*_{j_b}(t_b)\, \psi_{j_a}(t_a) \, v_0,v_0\right) \,,
& a > b \,.\\
\end{cases}
\]
Notice that $a\leq b$ iff $t_a\geq t_b$ and $a>b$ iff $t_a<t_b$. Now,
substitute (\ref{p}) and (\ref{p-star}) into $K_{ab}$ and take into
account (\ref{pmel}). This proves the formula for correlation
functions.

\end{proof}

\subsection{}Notice that operators $\psi(z,t)$ and $\psi^*(z,t)$
satisfy the difference equations:
\[
\psi(z,t+1)=(1-zx^+_{t+1/2})\psi(z,t), \ \ t\in D_+
\]
\[
\psi(z,t+1)=(1-z^{-1}x^-_{t+1/2})\psi(z,t), \ \ t\in D_-
\]
\[
\psi^*(z,t-1)=(1-zx^+_{t-1/2})\psi^*(z,t), \ \ t\in D_+
\]
\[
\psi^*(z,t-1)=(1-z^{-1}x^-_{t-1/2})\psi(z,t), \ \ t\in D_-
\]

These difference equations give the following difference equations for
correlation functions:
\begin{multline}\label{REC1}
K((t_1,h_1),(t_2,h_2))-K((t_1-1,h_1+1/2),(t_2,h_2)) \\
+x^+_{t_1-1/2}
K((t_1-1,h_1-1/2),(t_2,h_2))=\delta_{t_1,t_2}\delta_{h_1,h_2},
\  t_1\in D_+ \ ,
\end{multline}

\begin{multline}\label{REC2}
K((t_1,h_1),(t_2,h_2))-K((t_1-1,h_1-1/2),(t_2,h_2))\\
+x^-_{t_1-1/2}
K((t_1-1,h_1+1/2),(t_2,h_2))=\delta_{t_1,t_2}\delta_{h_1,h_2},
\ t_1\in D_- \ ,
\end{multline}

Using these equations and similar difference equations in $t_2$
one can express all correlation functions in terms of equal time
correlation functions.
\begin{remark}
The equations (\ref{REC1}) and (\ref{REC2}) together with
appropriate boundary conditions are the equations for the
inverse Kasteleyn matrix for the corresponding dimer model.
\end{remark}

\subsection{ The homogeneous restricted case}
In the homogeneous restricted case $0<q_t=q<1$ for $u_0\leq t\leq
u_N$ and $q_t=0$ otherwise.

In the homogeneous case $x^\pm_m=a^{\pm 1}q^{\pm m}$. The
partition function does not depend on $a$.  The functions
$\Phi_\pm(z,t)$ are:
\begin{eqnarray}
\Phi_+(z,t)&=&\prod_{m>t, m\in D^+}(1-zq^{m}a), \\
\Phi_-(z,t)&=&\prod_{m<t, m\in D^-}(1-z^{-1}q^{-m}a^{-1})
\end{eqnarray}

In this case we have
\[
R(t)=\left\{ \begin{array}{ll} a^{-1}q^{-v_i},& \ u_i<t<v_i\\
a^{-1}q^{-t},& \ v_{i-1}<t<u_i\end{array} \right.
\]
\[
R^*(t)=\left\{ \begin{array}{ll} a^{-1}q^{-t},& \ u_i<t<v_i\\
a^{-1}q^{-v_{i-1}},& \ v_{i-1}<t<u_i\end{array} \right.
\]

Notice, that when $q\to 0$ the density of horizontal tiles
converges to
\[
\rho(h,t)= \frac1{(2\pi i)^2} \int_{|z|=1+\epsilon}
\int_{|w|=1\-\epsilon} \frac1{z-w} \frac{dz \, dw}{z^{h+B(t)+\oh}
w^{-h-B(t)+\frac12}}\,.
\]
This integral is $1$ when $(h,t)$ is on a ``floor'' and is $0$
when it is on the ``wall''.

\section{The thermodynamic limit}\label{contlim}

Here we will study the limit $q\to 1$ of the homogeneous Gauss
distribution on restricted skew plane partitions when the number
of corners in the inner shape of diagrams remain finite.

We assume that $q=\exp(-r)$, \ $r\to +0$ and $U_i=ru_i, \
V_i=rv_i$ and  $ N$ remain fixed and that
$U_0<V_1<U_1<\dots<V_N<U_N$.

\subsection{Asymptotics of the partition function}

It is easy to compute the free energy of the system in this limit:
\begin{multline*}
F=-\log Z=- \sum_{ m<n, m \in D^-,
n\in D^+} \log(1-q^{n-m})= \\
-\frac{1}{r^2}\Bigg(\sum_{1\leq i\leq j\leq N}
  \int_{U_{i-1}<\mu<V_i}\int_{V_j<\nu<U_j}
\log (1-e^{\mu-\nu})d\mu d\nu\Bigg)+o\Big(\frac{1}{r^2}\Big)
\end{multline*}

Similarly, one can compute the asymptotic of the average volume of
a 3D partition:
\begin{multline*}
\lang|\pi|\rang= q\frac{\partial}{\partial q}Z= \sum_{m<n, m\in  D^-
n\in D^+}\frac{n-m}{1-q^{n-m}}=\\
\frac{1}{r^3}\Bigg(\sum_{i\leq j}
  \int_{U_{i-1}<\mu < V_i}\int_{V_j<\nu < U_j}
\frac{\nu-\mu}{1-e^{\mu-\nu}}d\mu d\nu\Bigg)+ o\Big(\frac{1}{r^3}\Big)
\end{multline*}

The first formula reflects essentially two dimensional nature of
the problem. The second formula implies that $r^{-1}$ is the
characteristic length of the system when $r\to 0$. 

\subsection{The function $S(z)$} Now let
us analyze the correlation functions (\ref{main-corr2}) in the 
limit $r\to +0$. Since $r^{-1}$ is a characteristic scale of the
system  in this limit we assume $\tau=t_ir, \ \chi=h_ir$ 
 remain finite. Depending on the value of $(\chi,\tau)$ we will either
keep the differences $t_i-t_j$ and \ $h_i-h_j$  finite, 
or we  will scale them as appropriate powers of 
$r$.

 When $r\to +0$ the
functions in the integral defining correlation functions behave as
$$
\frac{\Phi_-(z,t)}{\Phi_+(z,t)}z^{-h-B(t)-\frac12}=
\exp\left(\frac{S(z)}{r}\right)F(z)(1+O(r))
$$
where
\[
S(z)=\int_{\mu <\tau, \ \mu \in D_-} \log (1-z^{-1}e^\mu)d\mu
-\int_{\mu > \tau, \  \mu \in D_+} \log (1-ze^{-\mu})d\mu- (\chi
+B(\tau))\ln(z)
\]
and $F(z)$ can be computed explicitly.

In this limit the integral \eqref{main-corr2} becomes

\begin{multline}\label{Corr-as}
K(h_1,t_1),(h_2,t_2))=\frac{1}{(2\pi i)^2}\times  \\
\int_{C_z} \int_{C_w}
e^{\tfrac{S(z;\chi_1,\tau_1)-S(w,\chi_2,\tau_2)}{r}}\frac{F(z;\chi_1,\tau_1)}{F(w;\chi_2,\tau_2)}
\frac{\sqrt{zw}}{z-w}\frac{dz}{z}\frac{dw}{w}\, (1+o(1)) \,. 
\end{multline}

The integration contours are described in the previous section.
For example, if $N=2$ and $U_1<\tau_1,\tau_2<V_2$ the contours
are:
\[
C_z\times C_w=\left\{ \begin{array}{ll}
e^{V_2}> |z|, |w|>|z|, |w|>e^{\tau}& \mbox{ if $\tau_1>\tau_2$} \\
e^{V_2}>|z|>|w|>e^{\tau} & \mbox{ if $\tau_1\geq\tau_2$}
\end{array}\right.
\]

The integral (\ref{Corr-as})  can be computed by the steepest descent method.
In order to do this
one should first analyze critical points of $S(z)$ and then
deform contours of integration accordingly.

\section{Critical points of $S(z)$ and the deformation of contours}\label{crpoints}

The function $S(z)$ can be written as a sum of dilogarithms:
$$
S(z)=-(\chi +B(\tau))\ln(z)+\sum^N_{i=0} \Li_2(ze^{-U_i})
-\sum_{i=1}^N\Li_2(ze^{-V_i}) -\Li_2(ze^{-\tau}) \ .
$$
where
$$
\Li_2(z)=\int_0^z\frac{\ln(1-x)}{x}dx
$$

Critical points of $S(z)$ are zeros of
\[
z\frac{\partial}{\partial z}S(z)=-(\chi +B(\tau))+ \int_{\mu
> \tau, \  \mu \in D_+}\frac{ze^{-\mu}}{1-ze^{-\mu}}d\mu
+\int_{\mu  <\tau, \  \mu \in D_-}
\frac{z^{-1}e^{\mu}}{1-z^{-1}e^{\mu}}d\mu
\]

This is equivalent to the following equations.
\[
-(\chi +B(\tau))-\sum_{0\leq
 j<i}\log(\frac{1-z^{-1}e^{V_{j+1}}}{1-z^{-1}e^{U_j}})
+\log(\frac{1-ze^{-U_i}}{1-ze^{-\tau}})+\sum_{
 N\geq j>i}\log(\frac{1-ze^{-U_j}}{1-ze^{-V_j}})=0
\]
when $V_i<\tau<U_i$ and
\[
-(\chi +B(\tau))-\sum_{0\leq
 j<i}\log(\frac{1-z^{-1}e^{V_{j+1}}}{1-z^{-1}e^{U_j}})
-\log(\frac{1-z^{-1}e^{\tau}}{1-z^{-1}e^{U_i}})+\sum_{
 N\geq j>i}\log(\frac{1-ze^{-U_j}}{1-ze^{-V_j}})=0
\]
when $U_i<\tau <V_{i+1}$.

Exponentiating these equations we obtain
\begin{equation}\label{cp}
\exp(-\chi -B(\tau)-L(\tau)) \prod_{0\leq j\leq
N}(1-ze^{-U_j})=(1-ze^{-\tau})\prod_{1\leq
  j\leq N}(1-ze^{-V_j})
\end{equation}
where
\[
L(\tau)=\left\{\begin{array}{ll}
\sum_{0\leq j < i}(V_{j+1}-U_j) & \mbox{if $V_i< \tau <U_i$} \\
\sum_{0\leq j< i}(V_{j+1}-U_j) +\tau- U_i& \mbox{if $U_i< \tau
<V_{i+1}$}
\end{array}
\right.
\]
The number $L(\tau)$ is the total length of $(-)$ intervals which
are to the left of $\tau$.

It is easy to see that
\[
L(t)+B(t)=t/2 - u_0 \,.
\]
Therefore the equation (\ref{cp}) can be written as
\begin{equation}\label{cp2}
e^{\chi-\tau/2} z - e^{\chi+\tau/2} = f(z)\,
\end{equation}
where
$$
f(z) = e^{U_0} \frac{\prod_0^N (z e^{-U_i} -1)} {\prod_1^N (z
e^{-V_i} -1)}  \,.
$$

\subsection{The number of roots}
\begin{theorem} \label{realroots}The equation (\ref{cp2}) has either $N$ real
solutions or $N-2$ real solutions and two complex conjugate.
\end{theorem}
\begin{proof}
The function $f(z)$ has simple poles at $z=v_i$ and $z=\infty$ and
simple zeros at $z=u_i$. A sample graph of the function $f(z)$ is
plotted in Figure \ref{fofz}

\begin{figure}[htbp]
  \begin{center}
    {\scalebox{0.5}{\includegraphics{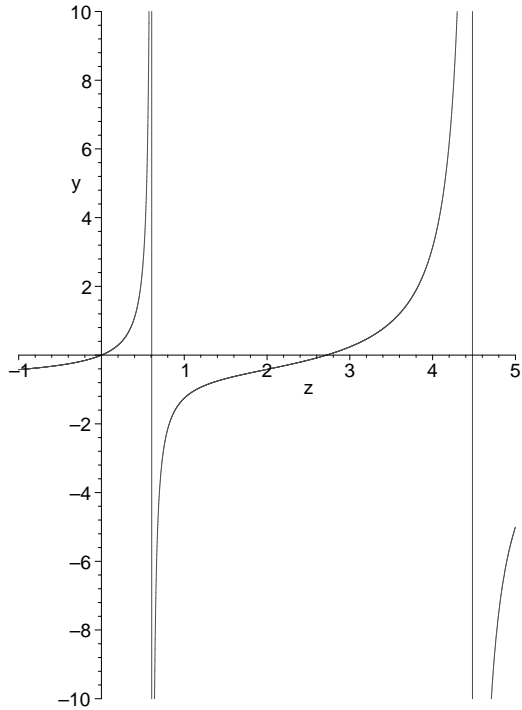}}}
    \caption{The graph of $f(z)$ for $\{U_0,V_1,U_1,V_2,U_2\}=\{-2,-0.5,1.5,2,3\}$}
    \label{fofz}
  \end{center}
\end{figure}

Solutions to the equation (\ref{cp2}) are intersection points of
the line $e^{\chi-\tau/}z-e^{\chi+\tau/2}$ with the graph of the
function $f(z)$.  Any straight line of non-infinite slope
obviously intersects the graph of $f(z)$ in at least $N-2$ points
(and in at most $N$ points, since the degree of $f$ equals $N$).
Hence among the roots of \eqref{cp2} there is at most one complex
conjugate pair.
\end{proof}

\subsection{The $N=1$ case} First, consider the equation
(\ref{cp}) When $N=1$ the normalization of $U$'s and $V$'s given
by the equation (\ref{norm-uv}) implies $V_1=0$. Introduce variables
\[
X=e^{\chi +\tau/2}, \ T=e^{-\tau},  \alpha=e^{U_0}, \
\beta=e^{-U_1},
\]
The equation for critical points of $S(z)$ is quadratic:
\begin{equation}\label{qe}
(\beta-XT) z^2+(X+XT-(1+\alpha\beta)) z + (\alpha-X)=0
\end{equation}
The discriminant of this equation is:
\begin{equation}\label{dis}
\Delta=
(X-XT-\alpha+\beta)^2-(X+XT)2(1-\alpha)(1-\beta)+(1-\alpha^2)(1-\beta^2).
\end{equation}

The structure of solutions depend on the value of the discriminant
$\Delta$.
\begin{itemize}
\item When $\Delta<0$ there are two complex conjugated critical
points

\item When $\Delta>0$ these critical points are real.

\item When $\Delta=0$ the two simple critical points degenerate
into one double critical point.

\end{itemize}

\subsection{The $N=2$ case}
This is the smallest value of $N$ when the function $S(z)$ can have
a triple critical point.

The function $S(z)$ for the case when $N=2$,  and $V_1+V_2<\tau<V_2$ is
\[
S(z)=\int_{U_0}^{v_1}\ln(1-z^{-1}e^\mu)d\mu+
\int_{v_1+v_2}^{\tau}\ln(1-z^{-1}e^\mu)d\mu-\int_{v_2}^{U_2}\ln(1-z^{-1}e^\mu)d\mu
-(\chi+B(\tau))\ln z
\]
We choose branches of logarithms such that the derivative of $S(z)$
has brunch cuts $[e^{U_0},e^{V_1}]$, \ $[e^{V_1+V_2}, e^\tau]$, and $[e^{V_2},
e^{U_2}]$.

The function $S(z)$ has three critical points. They are either all
real or there is a complex conjugate pair of simple complex critical points. 
Geometrically these critical points correspond to intersection points of
the line $e^{\chi-\tau/2}z-e^{\chi+\tau/2}$ with the graph of the
function $f(z)$. 

When the line $e^{\chi-\tau/}z-e^{\chi+\tau/2}$
intersects the graph of the function $f(z)$ transversally, the
intersection points are simple critical points of $S(z)$. 

At double critical points the the line
$e^{\chi-\tau/}z-e^{\chi+\tau/2}$ is tangent to $f(z)$ but it does  not
bisect the graph of $f(z)$ at the point where it is tangent to the graph. 

At a triple critical point the line is
tangent to the graph of $f(z)$ and bisects it.

By the definition, a double critical point $z$ of the function
$S(z)$ satisfies the two equations $S'(z)=0$ and $S''(z)=0$. We
have
\[
\zp S(z)=-\chi -\ln(ze^{-\frac{\tau}{2}}-e^{\frac{\tau}{2}})+\ln
f(z)
\]
\[
(\zp)^2S(z)=-\frac{ze^{-\frac{\tau}{2}}}{
ze^{-\frac{\tau}{2}}-e^{\frac{\tau}{2}}}+ z\frac{f'(z)}{f(z)}
\]
This gives the system of equations for double critical points:
\begin{eqnarray}\label{dcrit}
e^{\chi-\tau/2} z - e^{\chi+\tau/2}& = &f(z)\, \\
e^{\chi_0-\tau_0/2} & = & f'(z) \,
\end{eqnarray}

These equations define the curve in the $(\chi,\tau)$ plane. We
will say the point $(\chi,\tau)$ on this curve is {\it generic} if it is
not a triple critical point, i. e. if $S'''(z)\neq 0$ where $z$ is
the corresponding double critical point of $S(z)$ ( a solution to (\ref{dcrit}).

Denote by $z_0$ the triple critical point of $S(z)$ and by $(\chi_0,\tau_0)$
the corresponding values of $\chi$ and $\tau$. By definition
$S'(z_0)=S''(z_0)=S'''(z_0)=0$, which gives one more equation in
addition to (\ref{dcrit}) :
\[
f''(z_0)=0
\]

It is clear from the shape graph of $f(z)$ that for each value of
$\tau$ there are either $2$ or $4$ double critical points of
$S(z)$. They satisfy the following inequalities:
\begin{enumerate}
\item $-\infty <\tau < V_1$, then $z_1<0$ and $0<z_2<e^{V_1}$

\item $V_1<\tau <U_1$, then $z_1<0$ and $e^{V_1}<z_2<e^\tau$

\item $U_1<\tau <\tau_0$, then $z_1<0$, \ $e^{V_1}<z_2< z_0$,
 \ $z_0< z_3<e^{U_1}$, and $e^\tau <z_4< e^{V_2}$.

\item $\tau_0<\tau <V_2$, then $z_1<0$ and $e^\tau< z_2<e^{V_2}$.

\item $V_2<\tau< U_2$, then $z_1<0$ and $e^{V_2}<z_2$.

\end{enumerate}

It is also clear that if $\tau_0>U_1$, then
\[
e^{V_1}<z_0<e^{U_1}
\]

\subsection{Deformation of integration
contours}\label{cont-deform} We want to deform contours of
integration in (\ref{Corr-as}) to position them in the way the
steepest descent methods requires. If $z_c$ is a critical point of
$S(z)$, and if it does not lie on a branch cut of the function
$S(z)$, the integration contour should be deformed to a contour
which lies on the curve $\Im(S(z))=\Im(S(z_c))$.

The function $S(z)$ has branch cuts along the real line. However,
only $\exp(\frac{S(z)}{r})F(z)$, which is the leading term of the
asymptotic of $\frac{\Phi_-(z,t)}{\Phi_+(z,t)}z^{-h-B(t)-\frac12}$
appear in the integral (\ref{Corr-as}).

Zeros of $\Phi_+(z)$ are accumulating along $(e^{V_2},e^{U_2})$.
Therefore we can not deform $C_z$ through this segment but we can
deform through any other part of the real line.

Similarly, zeros of $\Phi_-(w)$ are accumulating in the segments
$(e^{U_0}, e^{V_1})$ and $(e^{U_1},e^\tau)$. Thus, the contour
$C_w$ can not be deformed through these segments but can be
deformed through any other segment of the real line.

\begin{figure}[htbp]
  \begin{center}
    \rotatebox{270}{\scalebox{0.3}{\includegraphics{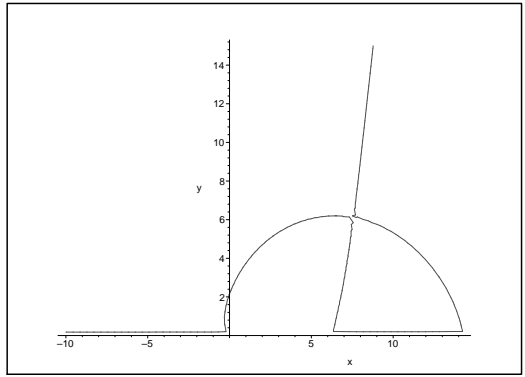}}}
    \caption{Level curves of $\Im(S)$ when all critical points are simple}
    \label{level-Im-simp}
  \end{center}
\end{figure}

\begin{figure}[htbp]
  \begin{center}
    \rotatebox{270}{\scalebox{0.3}{\includegraphics{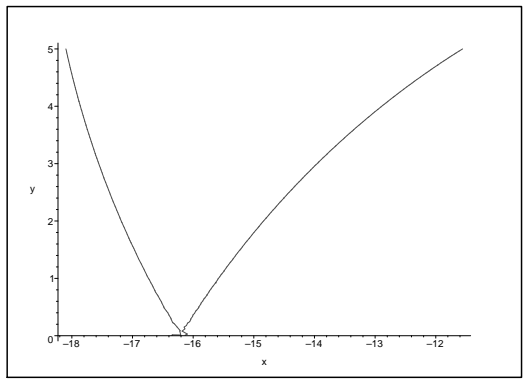}}}
    \caption{Level curves of $\Im(S)$ when there is a double critical
    point}
    \label{level-Im-double}
  \end{center}
\end{figure}

\begin{figure}[htbp]
  \begin{center}
    \rotatebox{270}{\scalebox{0.3}{\includegraphics{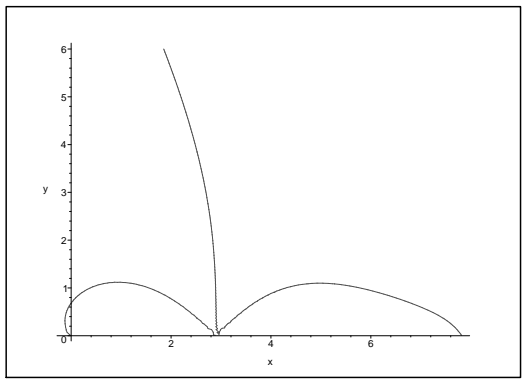}}}
    \caption{Level curves of $\Im(S)$ when there is a triple critical point}
    \label{level-Im-triple}
  \end{center}
\end{figure}

Therefore we have to deform contours $C_z$ and $C_w$ to the union
of appropriate branches of curves $\Im(S(z))=\Im(S(z_c+i0))$ and
$\Im(S(z))=\Im(S(z_c-i0))$. Figures \ref{level-Im-simp},
\ref{level-Im-double}, and \ref{level-Im-triple} show  these
curves for simple critical points, double critical points and the
triple critical point, respectively. 

\begin{figure}[htbp]
  \begin{center}
    {\scalebox{0.3}{\includegraphics{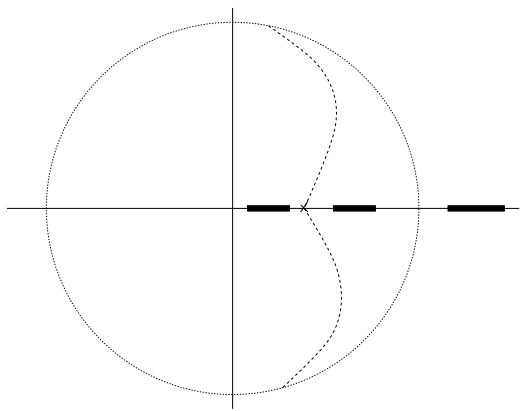}}}
    \caption{Deformation of the integration contour $C_z$ for
a double critical point}
    \label{cont-double-z}
  \end{center}
\end{figure}

\begin{figure}[htbp]
  \begin{center}
    {\scalebox{0.3}{\includegraphics{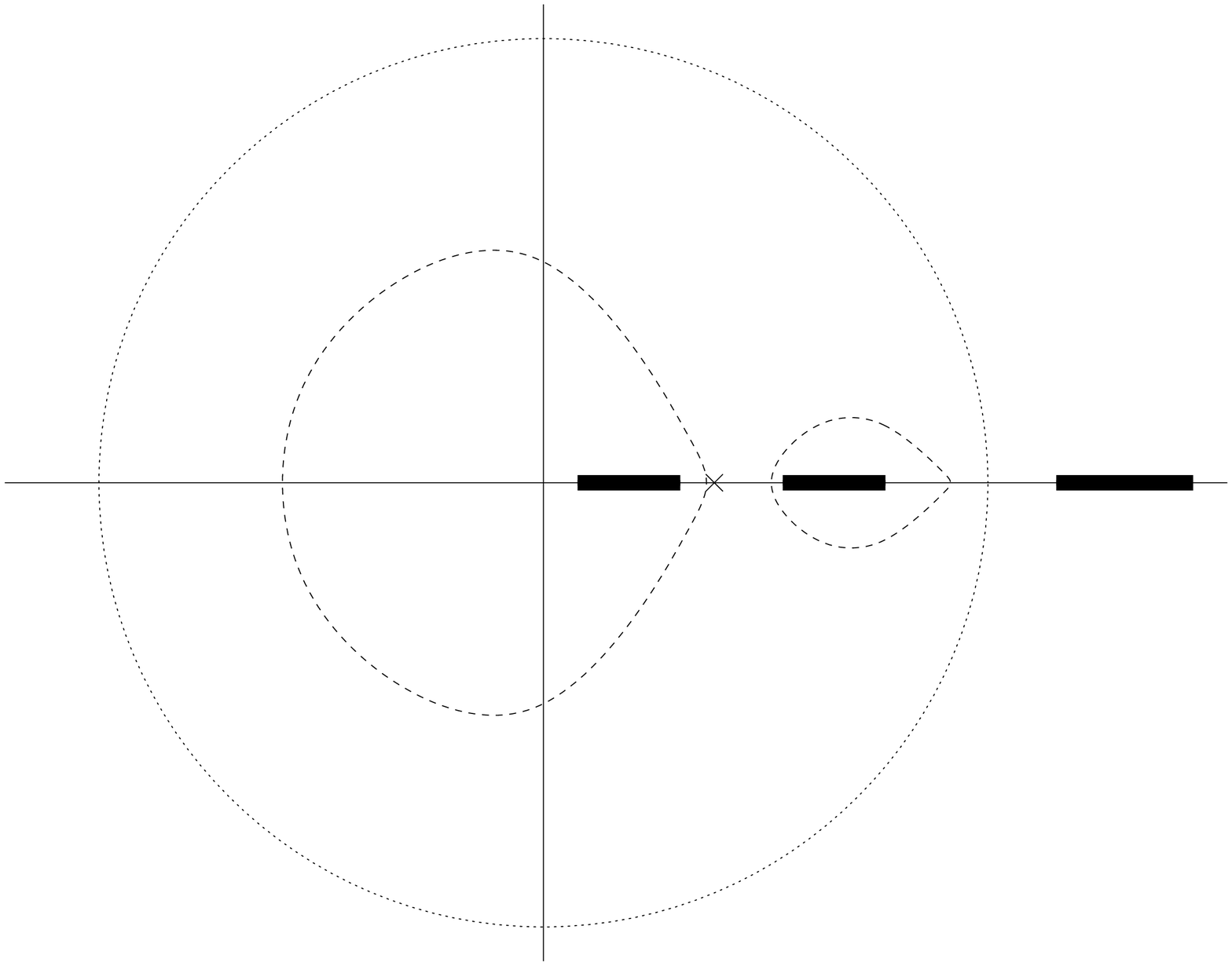}}}
   \caption{Deformation of the integration contour $C_w$ for
a double critical point}
    \label{cont-double-w}
  \end{center}
\end{figure}

\begin{figure}[htbp]
  \begin{center}
    {\scalebox{0.3}{\includegraphics{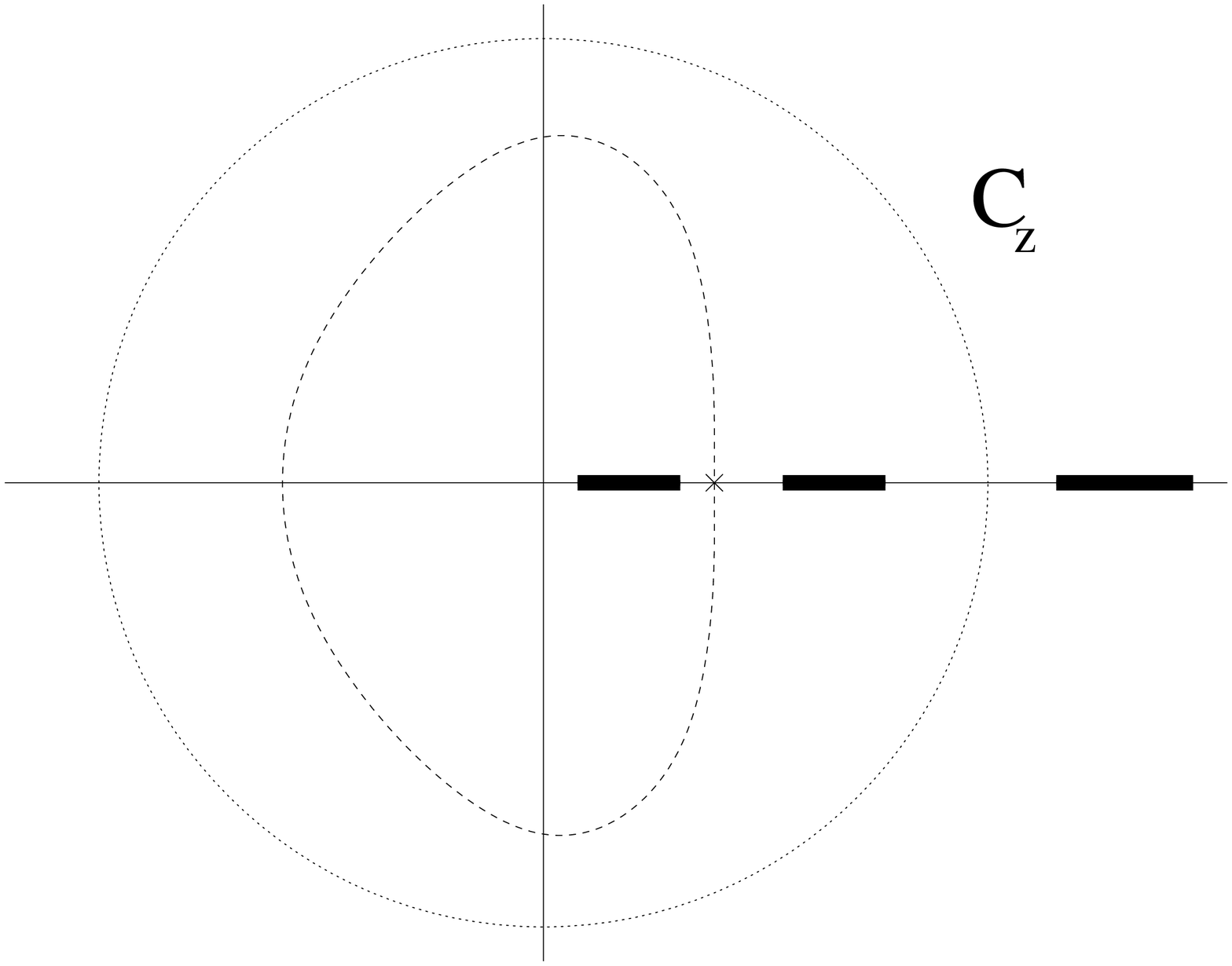}}}
    \caption{Deformation of the integration contour $C_z$ for a  triple
critical point}
    \label{cont-triple-z}
  \end{center}
\end{figure}

\begin{figure}[htbp]
  \begin{center}
    {\scalebox{0.3}{\includegraphics{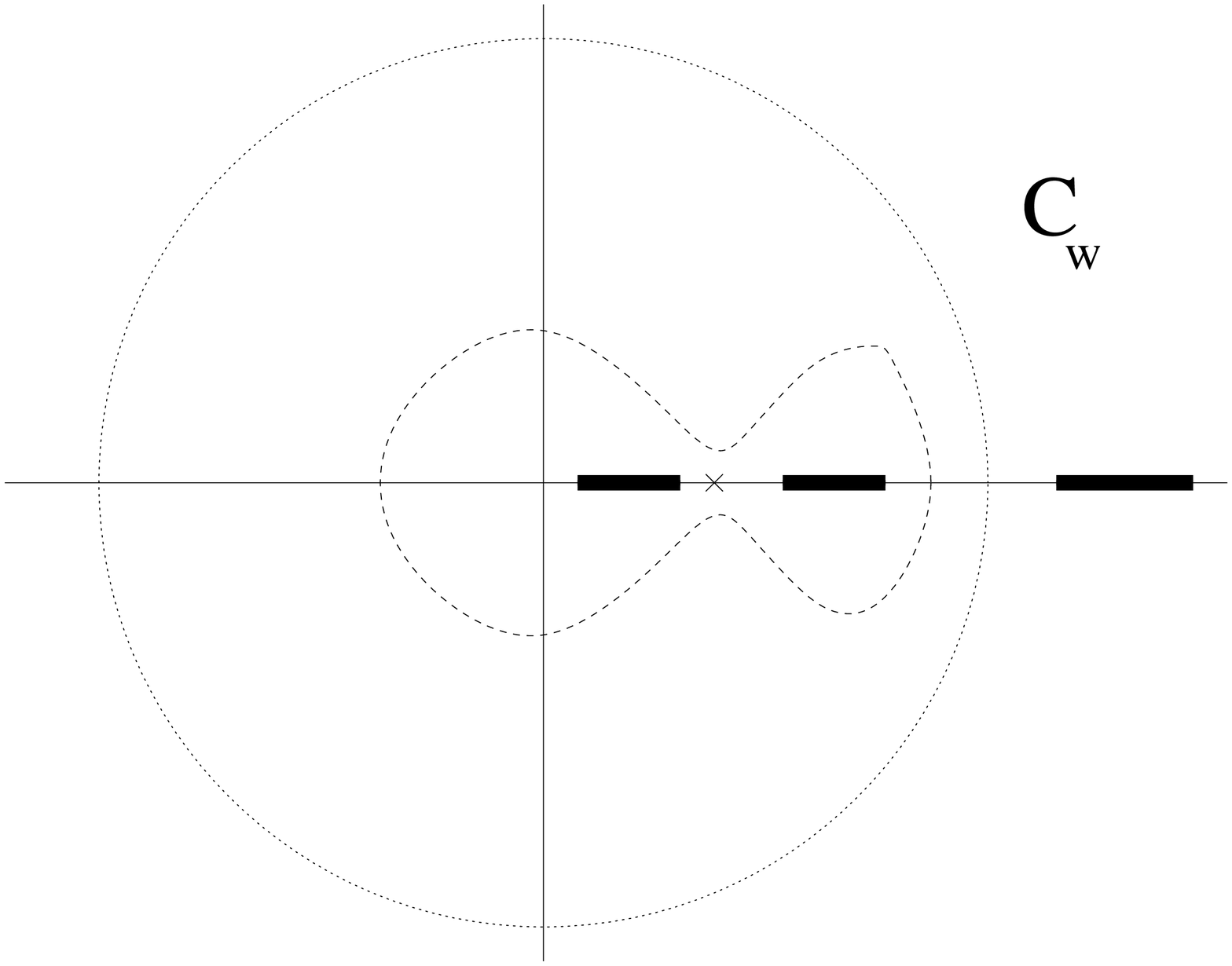}}}
    \caption{Deformation of the integration contour $C_w$ for a triple
    critical point}
    \label{cont-triple-w}
  \end{center}
\end{figure}

\begin{figure}[htbp]
  \begin{center}
    {\scalebox{0.3}{\includegraphics{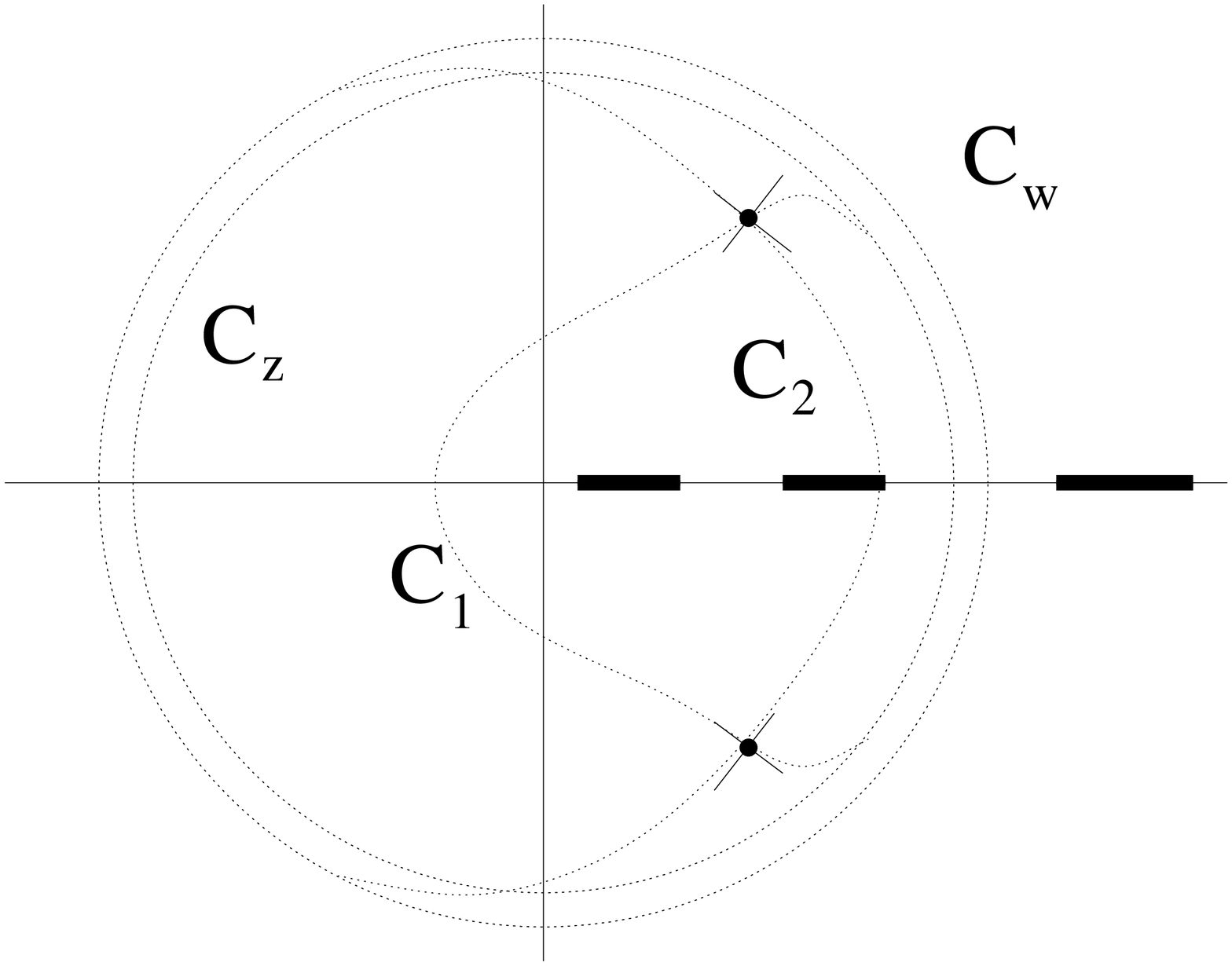}}}
   \caption{Deformation of the integration contour for
a simple critical point when $t_1\geq t_2$}
    \label{cont-simple-12}
  \end{center}
\end{figure}

\begin{figure}[htbp]
  \begin{center}
    {\scalebox{0.3}{\includegraphics{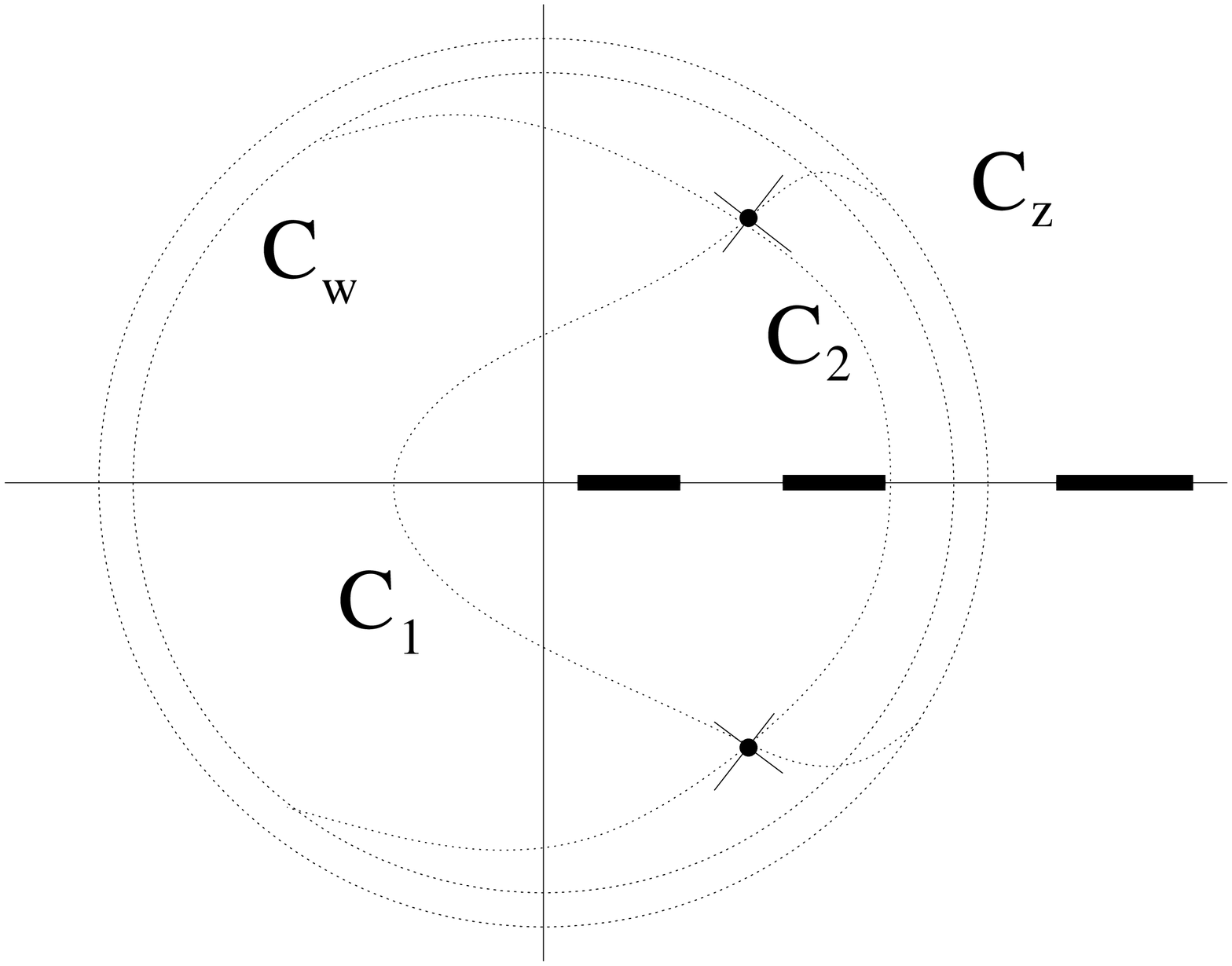}}}
   \caption{Deformation of the integration contours for
a simple critical point when $t_1< t_2$}
    \label{cont-simple-21}
  \end{center}
\end{figure}

Deformed contours of
integration $C_z$ and $C_w$ are shown in  Figures \ref{cont-simple-12},
\ref{cont-simple-21}, \ref{cont-double-z}, \ref{cont-double-w},
\ref{cont-triple-z}, \ref{cont-triple-w} for simple, double,
and triple critical points respectively.



While deforming contours $C_z$ and $C_w$ one should keep track on
the residues at the pole $z=w$. We will discuss this later.

\section{Simple critical points and bulk limit}

\subsection{Bulk limit of correlation functions}
Here we will compute the asymptotic of the integral
(\ref{Corr-as}) when there is a pair of complex conjugate critical
points and when $\Delta h=h_1-h_2$ and $\Delta t=t_1-t_2$ are
fixed in the limit $r\to 0$.

Deform contours $C_z$ and $C_w$ to $C_1$ and $C_2$ respectively as
it is  shown on Fig. \ref{cont-simple-12} for $t_1\geq t_2$ and 
as it is shown on Fig. \ref{cont-simple-21} for $t_1<t_2$.

Taking into account the residue at $z=w$ we have the identity:
\begin{equation}\label{simple-deform}
\int_{C_z}\int_{C_w}...=\int_{C_1}\int_{C_2}....+\int_{\gamma_\pm}...
\end{equation}
Here we take $+$ for $\Delta t\geq 0$ and $-$ otherwise. If
$U_i<\tau<V_{i+1}$ then $\gamma_+$ is a simple curve connecting
$z_c$ and $\bar{z_c}$ and passing through positive part of the
real line with $e^{V_{i+1}}<z$. Similarly, $\gamma_-$ is a simple
curve connecting $z_c$ and $\bar{z_c}$ and passing through the
negative part of the real line. If $V_i<\tau<U_i$ then
$\gamma_\pm$ are again simple curves connecting $z_c$ and
$\bar{z_c}$ in such a  way that $\gamma_-$ intersects the negative
part of the real line and $\gamma_+$ intersects the positive part
of the real line at $e^\tau<z$.

As $r\to 0$ only the second term in the right hand side of
(\ref{simple-deform}) will have finite limit. The first term will
vanish. The computations are similar to \cite{OR}. The integrals
along $\gamma_\pm$ converge as $r\to 0$ to
\begin{equation}
\kappa(\Delta h,\Delta
t)=\left\{\begin{array}{ll}B_{\gamma_+}^{(\e(\tau))}(\Delta
t,\Delta h+\e(\tau)\Delta t/2) &, \ \Delta t\geq 0
\\ B_{\gamma_-}^{(\e(\tau))}(\Delta t,\Delta h+\e(\tau)\Delta t/2) &, \ \Delta t<0
\end{array}\right.
\end{equation}\label{cf}
where
\[ B_\gamma^{(\pm)}(k,l)= \frac{1}{2\pi i}\int_{\gamma}
(1-e^{\mp\tau}w^{\pm 1})^{k} w^{-l-1} \, dw \ ,
\]
Here $\e(\tau)=1$ when $ \ V_i<\tau<U_i$ and $\e(\tau)=-1$ when
$U_i<\tau <V_{i+1}$. The contours $\gamma_\pm$ are as above.

For the limit of correlation functions we obtain the following
answer:
\[
\lim_{r\to 0}\lang\rho_{h_1,t_1}\dots \rho_{h_n,t_n}\rang=
\det(\kappa_{a,b})_{1\leq a,b\leq n}
\]
where
\[
\kappa_{a,b}=\kappa((h_a-h_b),(t_a-t_b))= \lim_{r\to +0}
K((t_a,h_a),(t_b,h_b))
\]

Correlation functions for finite $q$ satisfy recurrense equations
(\ref{REC1}) and (\ref{REC2}). In the limit these
recurrence relations turn into difference equations for correlation
functions (\ref{cf}):
\[
\kappa(h, t )-\kappa( h +\e(\tau) 1/2,  t-1)+
e^{-\e(\tau)\tau}\kappa(h -\e(\tau) 1/2 , t-1 )= \delta_{ h,
0}\delta_{ t, 0}
\]
where $\epsilon(\tau)=+1$ when $\tau\in D^+$, $\epsilon(\tau) =-1$
when $\tau \in D^-$. These equations can also be directly deduced
from the integral representation of correlation functions.

\begin{remark} The
difference operator in these equations is the Kasteleyn matrix
on the infinite hexagonal lattice. The pairwise correlation function
(\ref{cf}) can be regarded as the inverse for this matrix
with the boundary conditions determined by $(\tau, \chi)$.
\end{remark}

For one-time correlation functions we have:
\[
\kappa(h,0)=|z_c|^{-h}\frac{\sin(\theta h)}{\pi h}
\]
where $\theta$ is the argument  of $z_c$, that is,  $z_c=|z_c|e^{i\theta}$.

Notice that the factor $|z_c|^h$ does not contribute to the
equal-time density correlation functions and we have:

$$
\lim_{r\to 0}\lang\rho_{h_1,t}\dots \rho_{h_n,t}\rang=
\det\left(\frac{\sin(\theta(h_a-h_b))}{\pi (h_a-h_b)}
\right)_{1\leq a,b\leq
n}
$$
Here we assume that $h_1>h_2>\dots h_n$.

\subsection{Limit shape}
The form of the one-point correlation function implies that 
the density of horizontal tiles is
\[
\rho(\tau,\chi)=\frac{\theta}{\pi}
\]
where $\theta$ is the argument of $z_c$. The limit shape 
can reconstructed from this density by integration: 
\begin{equation}
z(\tau,\chi)=\int_{-b(\tau)}^\chi (1-\rho_(\tau,s)) \, ds \,.
\end{equation}
\begin{equation}
x(\tau,\chi)=z(\tau,\chi)-\chi-\frac\tau2\,, \quad
y(\tau,\chi)=z(\tau,\chi)-\chi+\frac\tau2\,,
\end{equation}
Here $b(\tau)=\max(-\tau/2+U_0, \tau/2+U_N)$.
Thus, the information about the limit shape is in the structure of
critical points of the function $S(z)$. 

When the point $(\tau,
\chi)$ is such that all critical points of $S(z)$ are real the
limit of correlation functions (\ref{cf}) is either $1$ or $0$. It
is zero if the maximum of $S(z)$ is inside of the cycle of
integration and it is $1$ if it is outside. The corresponding 
point $(x,y,z)$ lies on a facet (flat 
part of the limit shape). 

If $(\tau,\chi)$ is such that  there is a pair of complex
conjugate simple critical points, the point $(x,y,z)$
lies  in the disordered region (curved part of the limit 
shape). 

The frozen boundary (that is, the boundary between the 
disordered region and the facets)
limit shape correspond to $(\chi,\tau)$ for which 
there exists a real zero of $S(z)$ of the multiplicity at least
$2$. Cusps of the frozen boundary correspond to triple
critical points.

\begin{figure}[htbp]
  \begin{center}
    {\scalebox{0.3}{\includegraphics{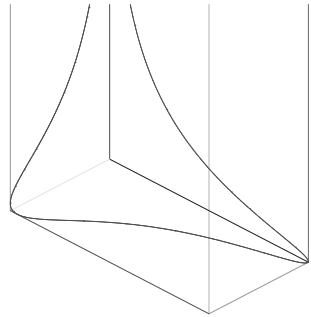}}}
   \caption{An example of frozen boundary for $N=1$}
    \label{bc1}
  \end{center}
\end{figure}

\begin{figure}[htbp]
  \begin{center}
    {\scalebox{0.3}{\rotatebox{-90}{\includegraphics{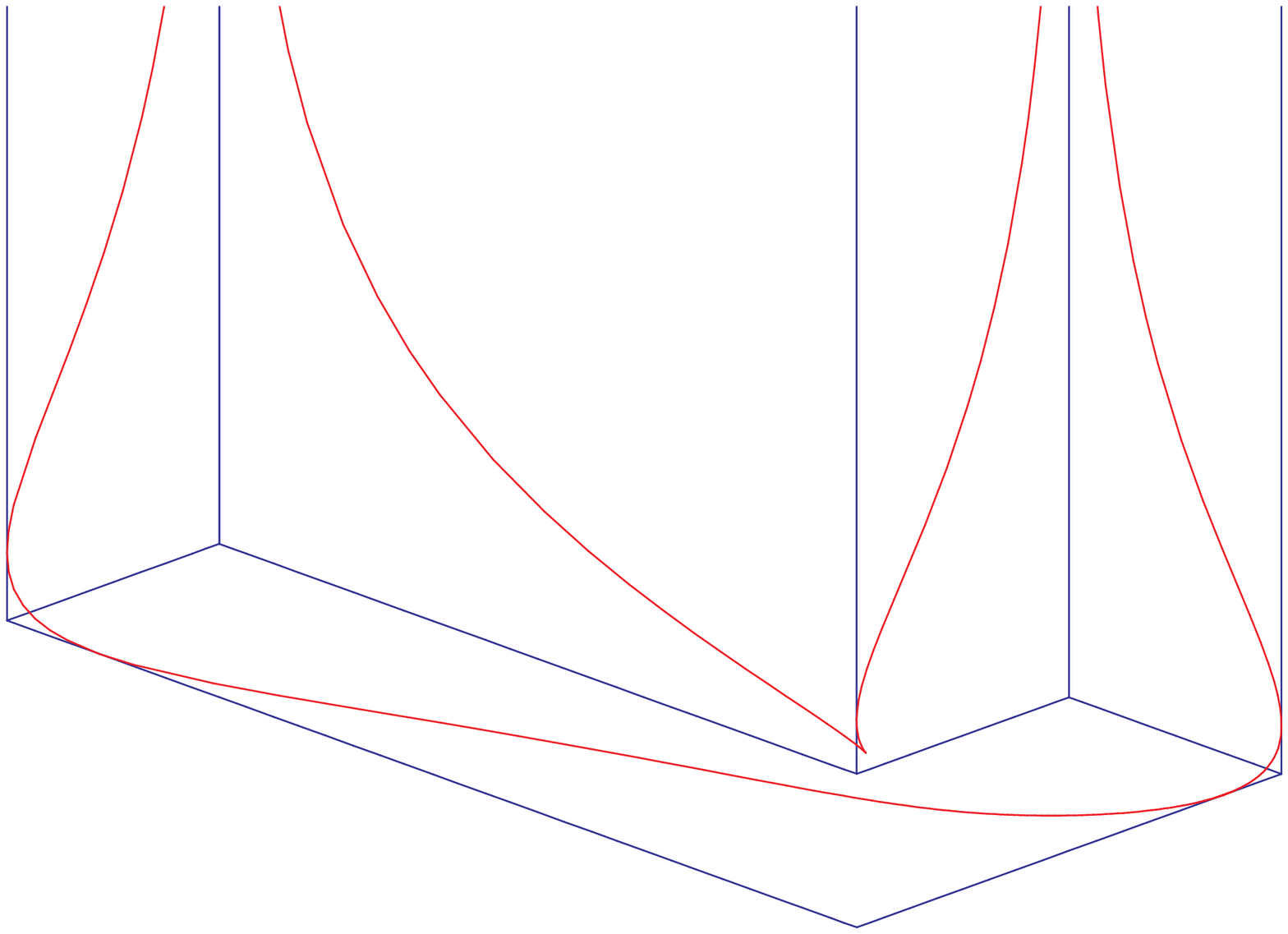}}}}
   \caption{An example of frozen boundary  for $N=2$}
    \label{shape1}
  \end{center}
\end{figure}

\begin{figure}[htbp]
  \begin{center}
    {\scalebox{0.3}{\rotatebox{-90}{\includegraphics{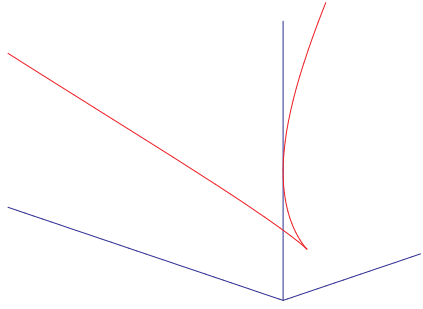}}}}
   \caption{The vicinity of a cusp}
    \label{shape2}
  \end{center}
\end{figure}

We plotted the frozen boundary in the 
$(\tau,\chi)$-plane on Fig. \ref{bc1} for $N=1$ and on
Fig. \ref{shape1}
for $N=1$. The vicinity of the cusp is magnified on Fig. \ref{shape2}.
It is instructive to compare these curves with the result of 
numeric simulation in Figure \ref{skew-random}.

\section{Double critical points and the scaling limit near the
boundary}
\subsection{}Now, let us assume that $(\chi_0,\tau_0)$ are such that
$z_0$ is a double critical point of $S(z)$. Consider the vicinity
of $(z_0,\chi_0, \tau_0)$ with coordinates $z=z_0\exp(\xi),
\chi=\chi_0+\delta\chi, \tau=\tau_0+\delta\tau$. Lowest degree
terms in the Taylor expansion around $(z_0,\chi_0,\tau_0)$ are:
\begin{eqnarray*}
S(z,\chi,\tau)&=&S(z_0,\chi_0,\tau_0)+\frac{A}{3!}\xi^3-(\delta\chi-B\delta\tau)\xi+
\frac{1}{2}C\delta\tau\xi^2 + \frac{1}{2}D\delta\tau^2\xi \\
&+& \frac{1}{3!}E\delta\tau^3
+H\delta\tau+\frac{1}{2}G\delta\tau^2-\ln(z_0)\delta\chi+ \dots
\end{eqnarray*}
where
\begin{equation}\label{B}
A= (\zp)^3S(z)|_{z_0,\tau_0}=\frac{z_0^2f''(z_0)}{f(z_0)}
\end{equation}
\begin{equation}\label{A}
B=z\frac{\partial^2}{\partial \tau\partial
z}S(z)|_{z_0,\tau_0}=
\frac{z_0e^{-\frac{\tau_0}{2}}+e^{\frac{\tau_0}{2}}}
{2(z_0e^{-\frac{\tau_0}{2}}-e^{\frac{\tau_0}{2}})}
\end{equation}
\begin{equation}\label{C}
C=(\zp)^2\frac{\partial}{\partial
\tau}S(z)|_{z_0,\tau_0}=-\frac{z_0}{(z_0e^{-\frac{\tau_0}{2}}-e^{\frac{\tau_0}{2}})^2}
\end{equation}
\begin{equation}\label{D}
D=z\frac{\partial}{\partial z}\frac{\partial^2}{\partial
\tau^2}S(z)|_{z_0,\tau_0}=\frac{z_0}{(z_0e^{-\frac{\tau_0}{2}}-e^{\frac{\tau_0}{2}})^2}
\end{equation}
\begin{equation}\label{E}
E= \frac{\partial^3}{\partial
\tau^3}S(z)|_{z_0,\tau_0}=-\frac{z_0}
{(z_0e^{-\frac{\tau_0}{2}}-e^{\frac{\tau_0}{2}})^2}
\end{equation}
\begin{equation}\label{H}
H=\frac{\partial}{\partial
\tau}S(z)|_{z_0,\tau_0}=\left\{\begin{array}{ll}
-\frac{1}{2}\ln(z_0)+\ln(z_0-e^{\tau_0)},& \ \tau_0\in D_- \\
-\frac{1}{2}\ln(z_0)+\ln(e^\tau-z_0)-\tau_0,& \ \tau_0\in D_+
\end{array}\right.
\end{equation}
\begin{equation}\label{G}
G=\frac{\partial^2}{\partial
\tau^2}S(z)|_{z_0,\tau_0}=\left\{\begin{array}{ll}
-\frac{e^{\tau_0}}{z_0-e^{\tau_0}},& \ \tau_0\in D_- \\
-\frac{z_0}{e^{\tau_0}-z_0},& \ \tau_0\in D_+
\end{array}\right.
\end{equation}
 Notice that $C=-D, E=-D, D>0, A/B>0$.
The sign of $A$ depends on the nature of the interface. It is
positive if the frozen region is above the melted region and it is
negative otherwise. It changes signs at the points $z_0=\exp(U_i)$
and at the triple critical points where $f''(z_0)=0$.

Rescaling local coordinates $\xi$, $\delta\chi$ and $\delta\tau$
as
\[
\xi=r^{\frac{1}{3}}\sigma, \ \
\delta\chi-B\delta\tau=r^{\frac{2}{3}}x, \ \
\delta\tau=r^{\frac{1}{3}}y \,
\]
we have
\[
\frac{S(z,\chi,\tau)-S(z_0,\chi_0,\tau_0)}{r}=\frac{A}{3!}\sigma^3-x\sigma+
\frac{1}{2}Cy\sigma^2+\frac{1}{2}Dy^2\sigma+\frac{1}{3!}Ey^3
+o(1)
\]

\subsection{}
As  $r\to 0$ the leading asymptotic of the integral
(\ref{Corr-as}) is determined by the leading asymptotic of the
integral:
\begin{equation}\label{Air-as}
K(h_1,t_1),(h_2,t_2))=\frac{1}{(2\pi i)^2}\int_{C_z} \int_{C_w}
e^{\frac{S(z;\chi_1,\tau_1)-S(w,\chi_2,\tau_2)}{r}}
\frac{\sqrt{zw}}{z-w}\frac{dz}{z}\frac{dw}{w}(1+O(r^{1/3})
\end{equation}
If $V_1<\tau<V_2$, the contour of integration is
\[
C_z\times C_w=\left\{ \begin{array}{ll}
1>|z|>|w|>e^{V_1} & \mbox{ if $V_1<\tau<U_1$} \\
e^{V_2}>|z|>|w|>e^{\tau} & \mbox{ if $U_1<\tau<V_2$}
\end{array}\right.
\]
for $t_1\geq t_2$. When $t_1<t_2$ the only difference is that
$|z|<|w|$.

Assume that coordinates $(h_i,t_i)$ are in the vicinity of the
point $(\chi_0,\tau_0)$, and that as $r\to 0$ they scale in the
following way:
\begin{equation}\label{Axyscale}
\chi_i=rh_i=\chi_0+r^{\frac{2}{3}}x_i-Br^{\frac{1}{3}}y_i, \
\tau_i=rt_i= \tau_0+r^{\frac{1}{3}}y_i \,.
\end{equation}
In this limit condition $t_1\geq t_2$, $t_2> t_1$ translate to
$y_1\geq y_2$ and $y_2>y_1$ respectively.

Deform contours of integration in such a way that they will pass
through this critical point and will follow the branches of curves
$\Im(S(z))=\Im(S(z_0))$ and $\Im(S(w))=\Im(S(z_0))$ where $\Re(S(z))$
and $-\Re(S(w))$ have local maximum. Deformed contours are described in
the section \ref{cont-deform}.

\begin{figure}[htbp]
  \begin{center}
    \scalebox{0.3}{\includegraphics{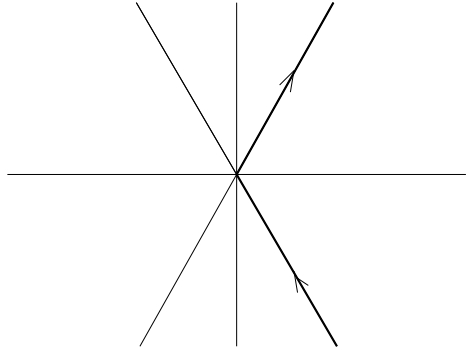}}
\caption{Integration contour for $\sigma$}
    \label{cont-airy-p.eps}
  \end{center}
\end{figure}

\begin{figure}[htbp]
  \begin{center}
    \scalebox{0.3}{\includegraphics{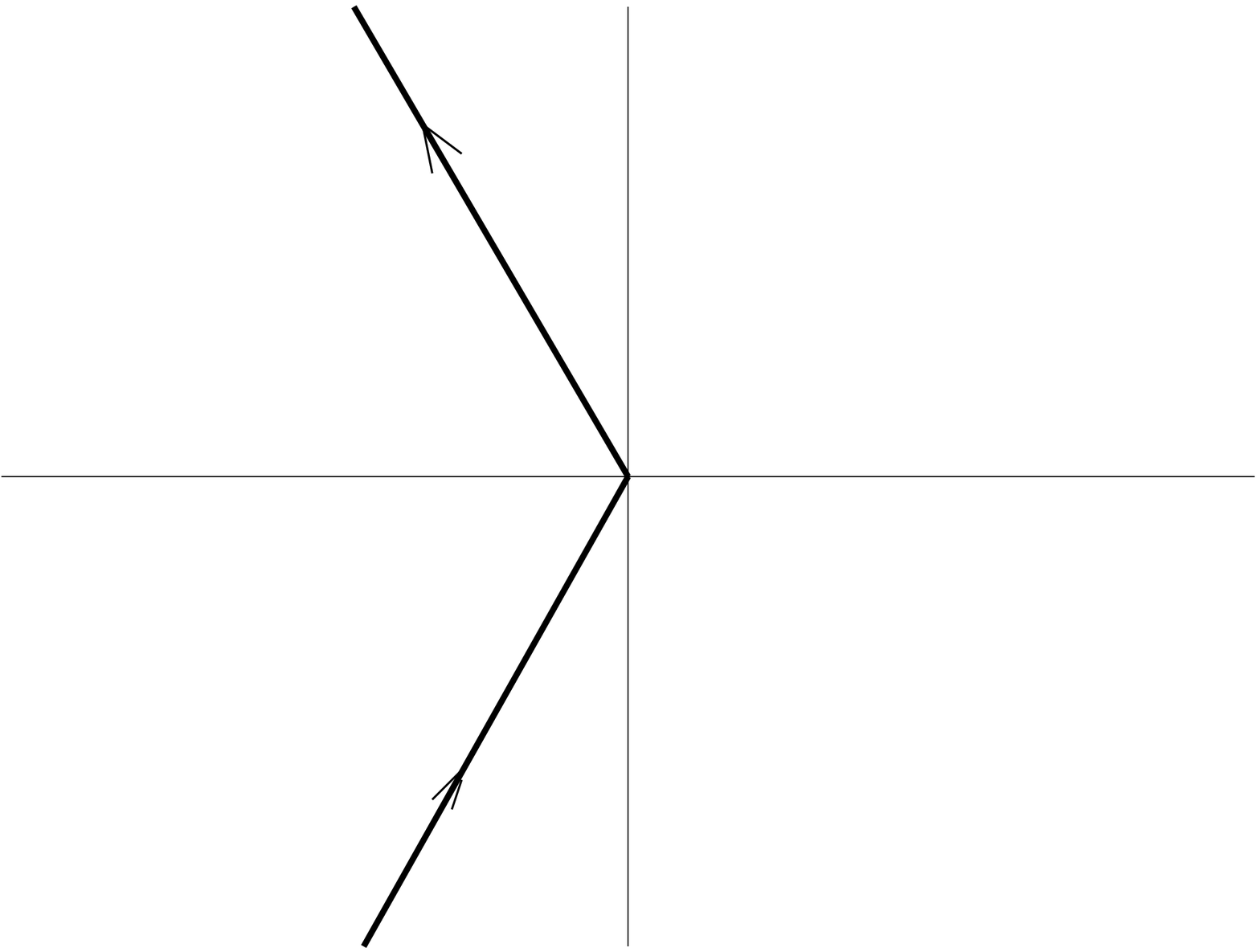}}
\caption{Integration contour for $\kappa$}
    \label{cont-airy-m.eps}
  \end{center}
\end{figure}

The leading contribution to
the asymptotic comes the vicinity of the point $z_0$ where
we will use local coordinates
$z=z_0(1+r^{\frac{1}{3}}\sigma)$ and
$w=z_0(1+r^{\frac{1}{3}}\kappa)$. The saddle point integration
contours for $\sigma$ and $\kappa$ in the limit $r\to 0$ are shown
on Fig. \ref{cont-airy-p.eps} and Fig. \ref{cont-airy-m.eps} for $A>0$.

The leading term of the asymptotic of the integral (\ref{Air-as})
is given by
\begin{multline}
K((h_1,t_1),(h_2,t_2))=
\exp(\frac{H}{r^{2/3}}(y_1-y_2)+\frac{G}{r^{1/3}}(y_1^2-y_2^2)
-\frac{\ln(z_0)}{r^{1/3}}(x_1-x_2)\\
-\frac{B\ln(z_0)}{r^{2/3}}(y_1-y_2)+\frac{1}{3!}E(y_1^3-y_2^3))
K_A((x_1-\frac{D}{2}y_1^2,y_1),(x_2-\frac{D}{2}y_2^2,y_2))(1+O(r^{1/3}))
\end{multline}
For $y_1\geq y_2$ the function $K_A$ is
\begin{multline}\label{K-Airy}
K_A((a_1,b_1),(a_2,b_2))=\\
\frac{1}{(2\pi i)^2}\int_{C_+}\int_{C_-}
\exp(\frac{A}{3!}(\sigma^3-\kappa^3)-\frac{D}{2}(b_1\sigma^2-b_2\kappa^2)-
a_1\sigma+a_2\kappa) \frac{d\sigma d\kappa}{\sigma-\kappa} \, .
\end{multline}
For $y_1<y_2$ it has an extra term which comes from the residue at
$\sigma=\tau$:
\begin{multline}\label{K-Airy1}
K_A((a_1,b_1),(a_2,b_2))= \frac{1}{(2\pi i)^2}\int_{C_+}\int_{C_-}
\exp(\frac{A}{3!}(\sigma^3-\kappa^3)-\frac{D}{2}(b_1\sigma^2-b_2\kappa^2)-\\
a_1\sigma+a_2\kappa) \frac{d\sigma d\kappa}{\sigma-\kappa} +
\frac{1}{2\pi i}\int_{C_+}\exp(-\frac{D}{2}(b_1-b_2)\sigma^2-
(a_1-a_2)\sigma) d\sigma  \, .
\end{multline}

The function $K_A((x_1,y_1),(x_2,y_2))$ is discontinuous at
$y_1=y_2$. The discontinuity is determined by the second term. It
is not difficult to see that it is equal to
\[
\sqrt{\frac{1}{2\pi
D|y_1-y_2|}}\exp(-\frac{(x_1-x_2)^2}{2D|y_1-y_2|})
\]

Notice that since correlation functions (\ref{corr}) are given by
determinants, the exponential factors in front of the integral in
(\ref{K-Airy}) will be canceled and we have:
\[
\lim_{r\to 0}\lang\rho_{h_1,t_1}\dots
\rho_{h_n,t_n}\rang=\det(K_A((x_a-\frac{D}{2}y_a^2,y_a),(x_b-\frac{D}{2}y_b^2,y_b,y_b)))
\]
Here we assume that $(h_i,t_i)$ are scaled as in \ref{Axyscale}.

It is easy to verify that
\[
K_{A}((x_1,0),(x_2,0))=\frac{\Ai(\lambda x_1)\Ai'(\lambda x_2)-
\Ai(\lambda x_2)\Ai'(\lambda x_1)}{\lambda(x_1-x_2)}
\]
where $\Ai(x)$ is the Airy function:
\[
\Ai(x)=\frac{1}{2\pi}\int_{\R}\exp(\frac{i}{3}t^3+ixt)dt ,
\]
and $\lambda=(\frac{2}{A})^{1/3}$.

For the scaling limit of the density of horizontal tiles
we have:
\begin{multline}
\rho_{A}(x,y)=\\-i\frac{1}{(2\pi)^2\lambda}\int_{\R}\int_{\R}
\exp\left(-i\frac{s^3-t^3}{3!}-\frac{i}{A^{1/3}}\left(\frac{C^2y^2}{2A}-x\right)(s-t)\right)
\frac{ ds dt}{s-t-i 0} \ .
\end{multline}

It can also be written as
\[
\rho_{A}(x,y)=\int_{-\infty}^{\lambda(\frac{C^2y^2}{2A}-x)}\Ai(v)^2dv
\]
In this form it is clear that the density function is positive.

\section{Triple critical points and the scaling limit near the
cusp}
\subsection{The singularity of the limit curve}
Then the following holds if $z_0$ is a triple critical point
\begin{equation}\label{s4}
(\zp)^4S(z)|_{z_o}=z_0^3\frac{f^{(3)}(z_0)}{f(z_0)}
\end{equation}
The proof is straightforward. First, one computes the fourth
derivative:
\[
(\zp)^4S(z)=\frac{z}{(ze^{-\frac{\tau}{2}}-e^{\frac{\tau}{2}})^2}-
\frac{2z^2e^{-\frac{\tau}{2}}}{(ze^{-\frac{\tau}{2}}-e^{\frac{\tau}{2}})^3}+
(z\frac{f'}{f}+3z^2(\frac{f'}{f})'+z^3(\frac{f'}{f})'')
\]
This expression reduces to (\ref{s4}) after taking into account
equations for $z_0$. The second derivative of $f(z)$ vanishes at
$z_0$ and, as it follows from the graph of $f(z)$, it is negative
for $z<z_0$ and positive for $z>z_0$. This implies the positivity
of $f^{(3)}(z_0)$. This means that the sign of $A$ is determined
by the sign of $f(z_0)$. The later is negative if $z_0<\exp(U_1)$
( when the cusp on the limit shape is turned to the right) and
is positive otherwise (i.e. when the cusp is turned to the left).

Now, let us find the behavior of the boundary curve near the cusp.
Assume that $z_0$  is a triple critical point corresponding to the
cusp with singularity at $(\chi_0,\tau_0)$, i.e. $z_0$,$a_0$ and
$b_0$ satisfy the system:
\begin{eqnarray}\label{c2}
e^{\chi_0-\tau_0/2}&=&f'(z_0), \\
e^{\chi_0+\tau_0/2}&=&z_0f'(z_0)-f(z_0), \\
0&=&f''(z_0)  .
\end{eqnarray}

If a point $(\chi,\tau)$ is at the boundary curve, it satisfies the
equations:
\begin{equation}\label{c1}
e^{\chi-\tau/2}=f'(z), \ \ e^{\chi+\tau/2}=zf'(z)-f(z),
\end{equation}
Let $(z,\chi,\tau)$ be a double critical point in a small vicinity of $(z_0,\chi_0,\tau_0)$:
\[
\chi=\chi_0+\delta\chi, \ \tau=\tau_0+\delta\tau, \
z=z_0+\e
\]
Then from the equations (\ref{c1},\ref{c2}) we obtain the
following asymptotic of the boundary curve near the cusp:
\[
\exp(\chi_0-\tau_0/2)(\delta\chi-\delta\tau/2)=\frac12f^{(3)}(z_0)\e^2+\frac{1}{6}f^{(4)}(z_0)\e^3+O(\e^4)
\]
\[
\exp(\chi_0+\tau_0/2)(\delta\chi+\delta\tau/2)=\frac12z_0f^{(3)}(z_0)\e^2+\frac{1}{6}(z_0f^{(4)}(z_0)+2f^{(3)}(z_0))\e^3+O(\e^4)
\]
when $\e\to 0$. From here we have:
\begin{multline}
\delta\chi=\frac{1}{4}\exp(-\chi_0-\frac{\tau_0}{2})((z_0+\exp(\tau_0))f^{(3)}(z_0)\e^2+\\
\frac{1}{3}((z_0+\exp(\tau_0))
f^{(4)}(z_0)+2f^{(3)}(z_0))\e^3+...)
\end{multline}
\begin{multline}
\delta\tau=\frac{1}{2}\exp(-\chi_0-\frac{\tau_0}{2})((z_0-\exp(\tau_0))f^{(3)}(z_0)\e^2+\\
\frac{1}{3}((z_0-\exp(\tau_0))
f^{(4)}(z_0)-2f^{(3)}(z_0))\e^3+...)
\end{multline}

After reparametrization
$$\e=\sigma-\frac{(z_0-e^{\tau_0})f^{(4)}(z_0)-2f^{(3)}(z_0)}
{6f^{(3)}(z_0)(z_0-e^{\tau_0})}$$ we have:
\[
\delta\tau=\frac{1}{2}\exp(-\chi_0-\frac{\tau_0}{2})(z_0-\exp(\tau_0))
f^{(3)}(z_0)\sigma^2+O(\sigma^3)
\ ,
\]
\[
\delta\chi=\frac{1}{4}\exp(-\chi_0-\frac{\tau_0}{2})((z_0+\exp(\tau_0))
f^{(3)}(z_0)\sigma^2-\frac{4}{3}\frac{\exp(\tau_0)}{z_0-\exp(\tau_0)}
f^{(3)}(z_0)\sigma^3)+O(\sigma^4) \ .
\]

This is a parametrization of a cusp singularity in the boundary of the 
limit shape.

\subsection{ The asymptotic of correlation functions}

Expanding the function $S(z,\chi,\tau)$ near the triple critical
point $z_0$ we obtain the following lowest degree terms of the
Taylor expansion:
\begin{multline*}
S(z,\chi_0+\delta\chi,\tau_0+\delta\tau)=S(z_0,\chi_0,\tau_0)+
\frac{A}{4!}\xi^4-\delta\chi\xi+B\delta\tau\xi+\frac{1}{2}C\delta\tau\xi^2+\\
+\frac{1}{2}D\delta\tau^2-\ln(z_0)\delta\chi+H\delta\tau+o(1)
\end{multline*}
where
\[
A=z_0^3\frac{f^{(3)}(z_0)}{f(z_0)}
\]
and $B$ and $C$ are as before, given by (\ref{B})(\ref{C}). Recall
that $A<0$ is the cusp in the limit shape is turned right and
$A>0$ if it is turned left.

Scaling local coordinates as
\[
\xi=r^{\frac{1}{4}}\sigma, \ \delta\chi -B\delta\tau
=r^{\frac{3}{4}}x, \ \delta\tau=r^\frac{1}{2}y \,
\]
we obtain the following asymptotic for $S$:
\begin{multline*}
\frac{S(z,\chi,\tau)-S(z_0,\chi_0,\tau_0)}{r}=\frac{A}{4!}\sigma^4
-x\sigma+ \frac{C}{2}y\sigma^2 +
\frac{D}{2}y^2+\\ \frac{H}{r^{1/2}}y-\frac{\ln(z_0)}{r^{1/4}}x+
\frac{B\ln(z_0)}{r^{3/4}}y+o(1)
\end{multline*}

Now let us find the asymptotic of the integral (\ref{Air-as}) as
$r\to 0$ assuming that coordinates are scaled as
\begin{equation}\label{Pxyscale}
rh_i=\chi_0+r^{\frac{3}{4}}x_i-Br^\frac{1}{4}y_i, \
rt_i=\tau_0+r^\frac{1}{2}y_i \, .
\end{equation}
where $(\chi_0, \tau_0)$ are coordinates of the tip of the cusp.

The asymptotic of the integral (\ref{Air-as}) as $r\to 0$ and 
\[
rh_i=\chi_0+r^{\frac{3}{4}}x_i-Br^\frac{1}{4}y_i, \
rt_i=\tau_0+r^\frac{1}{2}y_i \, .
\]
can be evaluated by the steepest descent method.
It is determined by the contribution
from the triple critical point $z_0$ and after deforming contours
of integration as it described in section \ref{cont-deform} the
leading term of the asymptotic is given by the integral
\begin{multline}
K((h_1,t_1)(h_2,t_2))=
\exp(-\frac{\ln(z_0)}{r^{1/4}}(x_1-x_2)+\frac{H}{r^{1/2}}(y_1-y_2)+\frac{B\ln(z_0)}{r^{3/4}}(y_1-y_2)+\\
\frac{D}{2}(y_1^2-y_2^2))K_P((x_1,y_1),(x_2,y_2)) (1+O(r^{1/4})
\end{multline}
where for $y_1>y_2$ we have
\begin{multline}
K_P((x_1,y_1),(x_2,y_2))=\\
\frac{1}{(2\pi i)^2}\int_{C_1}\int_{C_2}
\exp(\frac{A}{4!}(\sigma^4-\kappa^4)
+\frac{C}{2}(y_1\sigma^2-y_2\kappa^2)-x_1\sigma+x_2\kappa)
\frac{d\sigma d\kappa}{\sigma -\kappa} +\\
\sqrt{\frac{1}{2\pi
D|y_1-y_2|}}\exp(-\frac{(x_1-x_2)^2}{2D|y_1-y_2|})\, .
\end{multline}

For $y_1<y_2$ there is an extra term coming from the residue at
$\sigma=\tau$:

\begin{multline}
K_P((x_1,y_1),(x_2,y_2))=\\
\frac{1}{(2\pi i)^2}\int_{C_1}\int_{C_2}
\exp(\frac{A}{4!}(\sigma^4-\kappa^4)
+\frac{C}{2}(y_1\sigma^2-y_2\kappa^2)-x_1\sigma+x_2\kappa)
\frac{d\sigma d\kappa}{\sigma -\kappa} \, .
\end{multline}

\begin{figure}[htbp]
  \begin{center}
    \scalebox{0.3}{\includegraphics{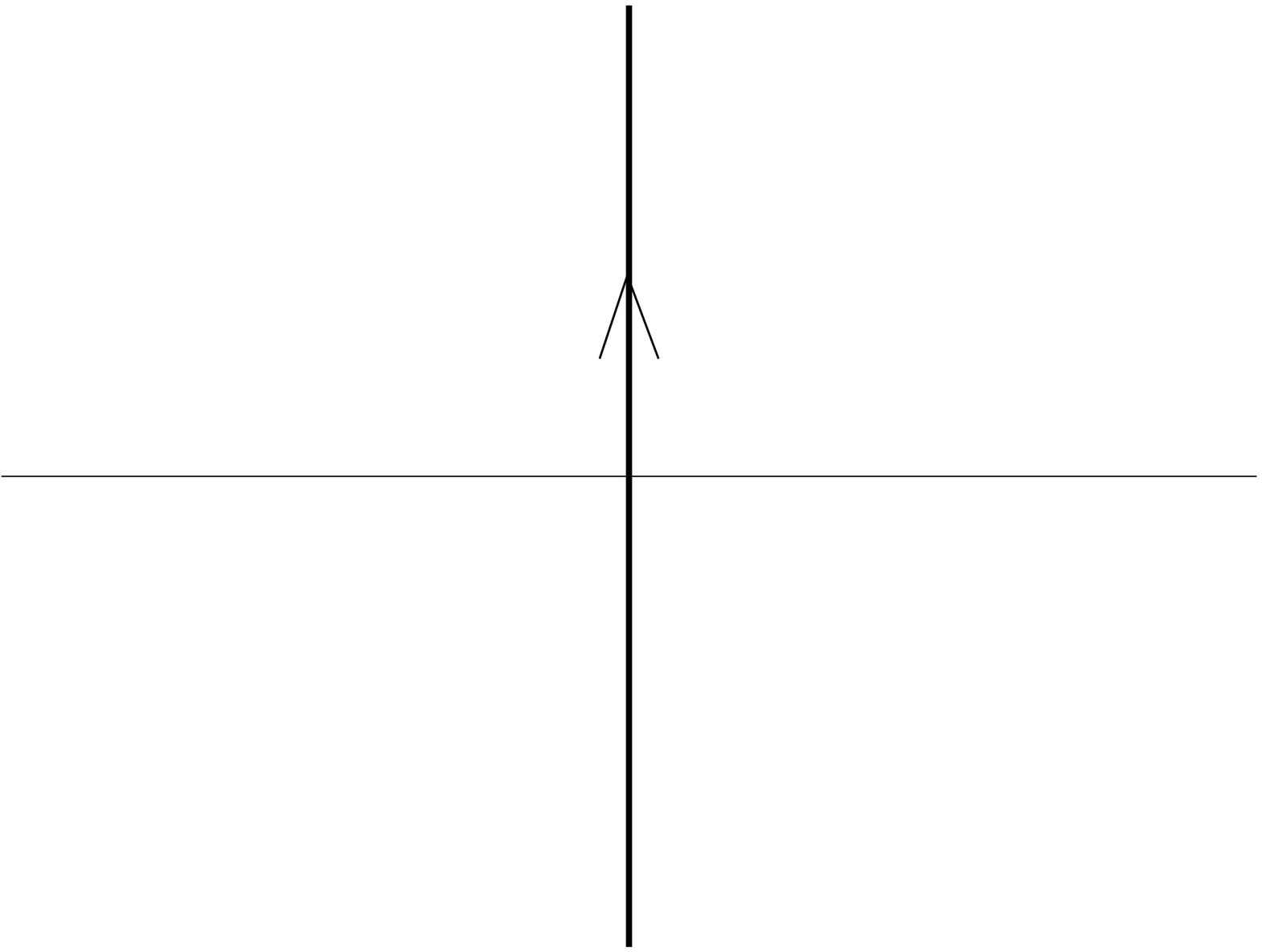}}
    \caption{Integration contour  for $\kappa$}
    \label{cont-per-c1.eps}
  \end{center}
\end{figure}

\begin{figure}[htbp]
  \begin{center}
    \scalebox{0.3}{\includegraphics{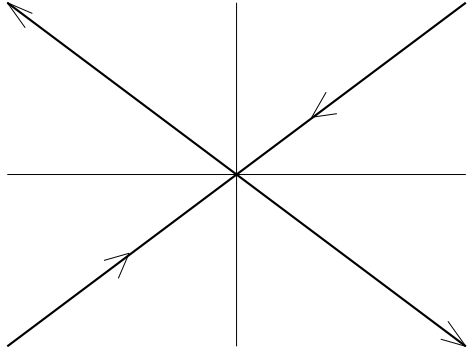}}
\caption{Integration contour for $\sigma$}
    \label{cont-per-c2.eps}
  \end{center}
\end{figure}

Integration contours for $A>0$ are shown on Fig.
\ref{cont-per-c2.eps}, and 
\ref{cont-per-c1.eps}.

Notice that the exponential factors cancels in the limit of
correlation functions and we have:

\[
\lim_{r\to 0}\lang\rho_{h_1,t_1}\dots
\rho_{h_n,t_n}\rang=\det(K_P((x_a,y_a),(x_b,y_b)))
\]
Here we assume $(h_i,t_i)$ scale  as in \ref{Pxyscale}).

\section{Some properties of Pearcey kernel}

Recall that the asymptotic of correlation functions is determined by the
following integral, which we will call Pearcey kernel

\begin{multline}
K_P((x_1,y_1),(x_2,y_2))=\\
\frac{1}{(2\pi i)^2}\int_{C_1}\int_{C_2}
\exp(\frac{A}{4!}(\sigma^4-\kappa^4)
+\frac{C}{2}(y_1\sigma^2-y_2\kappa^2)-x_1\sigma+x_2\kappa)
\frac{d\sigma d\kappa}{\sigma -\kappa} \, .
\end{multline}

After the appropriate scaling of variables we can set $A=6$ and  $C=1$.
Below we will review some of its properties.

\subsection{}

The function $K_P( (x_1,y_1),(x_2,y_2)))$  (with $A=6$ and $C=1$)satisfies the following
differential  equations:
\[
(\frac{\partial}{\partial y_i}-\frac{1}{2}\frac{\partial^2}{{\partial x_i}^2})K_P=0
\]
for $i=1,2$ and 
\[
(\frac{\partial}{\partial x_1}+\frac{\partial}{\partial x_2})K_P=-P_+(x_1,y_1)P_-(x_2,y_2)
\]
where 
\[
P_\pm(x,y)=\frac{1}{2\pi i}\int_{C_\pm} \exp(\pm(\frac{\sigma^4}{4}+\frac{1}{2}y\sigma^2-x\sigma))d\sigma
\]

For functions $P_\pm$ we have:
\[
(\frac{\partial^3}{{\partial x}^3}\pm x+y\frac{\partial}{\partial x})P_\pm=0
\]
\[
(\pm \frac{\partial}{\partial y}-\frac{1}{2}\frac{\partial^2}{{\partial x^2}})P_\pm=0
\]

\subsection{}

The function $P_+(x,y)$ can be written as
\[
P_+(x,y)=\frac{1}{2\pi i}(-\int_{\omega \R_+}+\int_{\omega^{-1}\R_+}-\int_{\omega^{-3}\R_+}+
\int_{\omega^{3}\R_+})
\]

or, as
\[
P_+(x,y)=\frac{1}{\pi}\Im(\omega\int_o^\infty \exp(-\frac{s^4}{4}-\frac{i}{2}ys^2)(e^{x\omega s}-e^{-x\omega s})
ds)
\]

From this integral representation one can find a power series formula
for the integral:
\[
P_+(x,y)=\frac{2}{\pi}\sum_{n\geq o, m\geq o, n+m=0(2)}\frac{(-1)^{\frac{n+m}{2}}2^{n-1}\Gamma(\frac{n+m+1}{2})}
{(2n+1)!m!}x^{2n+1}(-y)^m
\]
or, 
\begin{eqnarray*}
P_+(x,y)=&\frac{1}{\pi}&\sum_{k,l\geq 0}\frac{(-1)^{k+l}2^{2k}\Gamma(k+l+1/2)}{(4k+1)!(2l)!}x^{4k+1}y^{2l}
+\\
&\frac{1}{\pi}&\sum_{k, l\geq 0}\frac{(-1)^{k+l}2^{2k+1}\Gamma(k+l+3/2)}{(4k+3)!(2l+1)!}x^{4k+3}y^{2l+1}
\end{eqnarray*}

Now, using the identity $\Gamma(k+1/2)=2^{-2k+1}\sqrt{\pi}\frac{(2k)!}{k!}$ we arrive to
the formula
\begin{eqnarray*}
P_+(x,y)=&\frac{2}{\sqrt{\pi}}&\sum_{k,l\geq 0}\frac{(-1)^{k+l}2^{-2l+1}(2k+2l)!}{(4k+1)!(2l)!(k+l)!}x^{4k+1}y^{2l}
+\\
&\frac{2}{\sqrt{\pi}}&\sum_{k,l\geq 0}\frac{(-1)^{k+l}2^{-2l}(2k+2l+2)!}{(4k+3)!(2l+1)!(k+l+1)!}x^{4k+3}y^{2l+1}
\end{eqnarray*}

\subsection{}

The integral
\[
P_-(x,y)=\frac{1}{2\pi i}\int_{-\infty}^\infty\exp(-\frac{t^4}{4}+\frac{1}{2}yt^2+ixt)dt
\]
can be easily expanded into a power series in $x$ and $y$:
\[
P_-(x,y)=\frac{1}{\pi}\sum_{n,m\geq 0}\frac{(-1)^n2^{n-3/2}\Gamma(\frac{n+m}{2}+1/4)}{(2n)!m!}x^{2n}y^m,
\]

\subsection{}
One-time correlation functions can be expressed in terms of $P_\pm$ as follows::
\begin{multline}\label{1-time}
K(x_1,x_2|y)= \\
\frac{P''_+(x_1,y)P_-(x_2,y)-P'_+(x_1,y)P'_-(x_2,y)+P_+(x_1,y)P''_-(x_2,y)+yP_+(x_1,y)P_-(x_2,y)}{x_1-x_2}
\end{multline}

Here $P'(x,y)=\frac{\partial}{\partial x}P(x,y)$. This identity follows from
\begin{multline}
\int_{C_1}\int_{C_2}(\frac{\partial}{\partial \sigma}+\frac{\partial}{\partial \kappa})\\
\exp(\frac{A}{4!}(\sigma^4-\kappa^4)
+\frac{C}{2}(y_1\sigma^2-y_2\kappa^2)-x_1\sigma+x_2\kappa)
\frac{d\sigma d\kappa}{\sigma -\kappa}=0 \, .
\end{multline}

Taking the limit $x_1\to x_2=x$ in (\ref{1-time}) we obtain the
following expression for the 
scaling limit of the density of tiles near the cusp:
\[
\rho(x,y)=P'_+(x,y)P''_-(x,y)-P''_+(x,y)P'_-(x,y)-xP_+(x,y)P_-(x,y)
\]

\begin{figure}[htbp]
  \begin{center}
    \scalebox{0.5}{\includegraphics{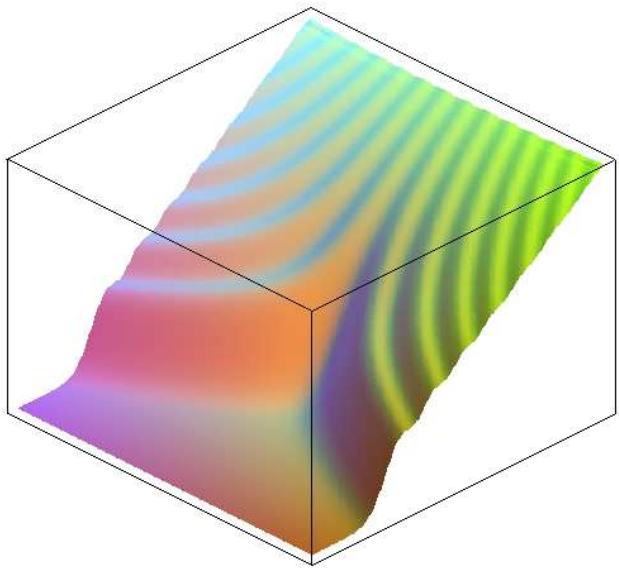}}
\caption{The scaling limit of the density of horizontal tiles near the
cups}
    \label{Density}
  \end{center}
\end{figure}

This function is plotted on Fig \ref{Density}.

\appendix

\section{Schur functions}\label{app_schur} 
\begin{figure}[htbp]
  \begin{center}
    \scalebox{0.3}{\includegraphics{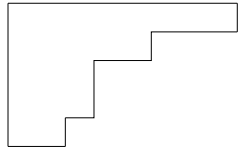}}
    \caption{A Young diagram}
    \label{yd.eps}
  \end{center}
\end{figure}

Recall that a partition is a
sequence of integers
\[
\lambda=(\lambda_1\ge\lambda_2\ge\dots\ge \lambda_n\ge 0)
\]
A diagram of a partition $\lambda$ (also known as its 
\emph{Young diagram}) has
$\lambda_1$ boxes in the upper row, $\lambda_2$ boxes in the next
row etc.,  see Figure 
\ref{yd.eps}. Unless this leads to a confusion, 
we identify partitions with their diagrams. The sum
\[
|\lambda|=\sum_{i\ge 1} \lambda_i
\]
is the area of the diagram $\lambda$.


A skew diagram $\lambda/\mu$ is a pair of two partitions
$\lambda$ and $\mu$ such that $\lambda_1\ge \mu_1, \
\lambda_2\ge \mu_2, \ \dots $. 
Graphically, it is obtained from $\lambda$ by removing $\mu_1$
first boxes from the first row, $\mu_2$ first boxes from the
second row etc. The size of the skew diagram $\lambda/\mu$ is
\[
|\lambda/\mu|=|\lambda|-|\mu|
\]

\begin{figure}[htbp]
  \begin{center}
    \scalebox{0.3}{\includegraphics{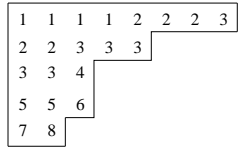}}
    \caption{A semi-standard  tableau}
    \label{sst.eps}
  \end{center}
\end{figure}

A semistandard tableau of the shape $\lambda/\mu$ with entries
$1,2,3,\dots$ is the result of writing numbers $1,2,3,\dots$ in
boxes of the diagram, one in each box, in such a way that the
numbers weakly descrease along the rows and strictly descrease 
along the colums, see an example in Figure \ref{sst.eps}


The Schur function corresponding to the skew  tableau
$\lambda/\mu$ is a symmetric function of variables $x_1,x_2,\dots$
which can be defined as the sum over all semi-standard tableaux of
shape $\lambda$ of monomials in $x_i$:
\[
s_{\lambda/\mu}(x_1,x_2,\dots)=\sum_{T_{\lambda/\mu}}
x_1^{m_1}x_2^{m_2}\dots
\]
where $m_i$ is the number of occurrences of $i$ in the tableau.

A semi-standard tableaux of shape $\lambda/\mu$ can be identified
with sequences of skew  diagrams in the following way.

Let us say that $\lambda\succ\mu$ if $\lambda$ and $\mu$
interlace, that is,
$$
\lambda_1 \ge \mu_1 \ge \lambda_2 \ge \mu_2 \ge \lambda_3 \ge
\dots \,.
$$
Let us call the sequence of  diagrams
$\lambda(1),\lambda(2),\dots \lambda(n)$ interlacing if
$\lambda(1)\succ\lambda_2\succ\dots\succ\lambda(n)$. It is clear
that if we will associate with such sequence the semi-standard
 tableau with $|\lambda(1)|-|\lambda(2)|$ entrees of $n$,
$|\lambda(2)|-|\lambda(3)|$ entries of $n-1$, etc. It is also
clear that this correspondence gives a bijection between
semi-standard tableaux of shape $\lambda(1)/\lambda(n)$ and
interlacing sequences of  diagrams which start with
$\lambda(1)$ and ends with $\lambda(n)$.

For skew Schur functions this bijection gives the following
formula
\begin{multline}\label{ssmon}
s_{\lambda/\mu}(x_1,\dots,x_n,0,0,\dots)= \\
\sum_{\lambda\succ\lambda(1)\succ\dots \succ\lambda(n-1)\succ\mu}
x_1^{|\lambda|-|\lambda(1)|}x_2^{|\lambda(1)| -|\lambda(2)|}\dots
x_n^{|\lambda(n-1)|-|\mu|}
\end{multline}
from now on we will write $s_{\lambda/\mu(x_1,x_2,\dots,x_n)}$ for
$s_\lambda(x_1,\dots,x_n,0,0,\dots)$.

Notice that
\begin{equation}\label{simpleschur}
s_{\lambda/\mu}(x_1,0,\dots)=\begin{cases}
x_1^{|\lambda|-|\mu|}\,, &  \mu \prec \lambda\,,\\
0\,, & \mu \not\prec \lambda\,.
\end{cases}
\end{equation}
and therefore we can write
\begin{multline}\label{recurschur}
s_{\lambda/\mu}(x_1,\dots,x_n)=\\
\sum_{\lambda\succ\lambda(1)\succ\dots \succ\lambda(n-1)\succ\mu}
s_{\lambda/\lambda(1)}(x_1) s_{\lambda(1)/\lambda(2)}(x_2)\dots
s_{\lambda(n-1)/\mu}(x_n)
\end{multline}

The function $s_\lambda(x_1,x_2,\dots,x_n)$ is the character of
the irreducible representation of $GL_n$ with the highest weight
$\lambda$ computed on the diagonal element with entries
$x_1,x_2,\dots,x_n$. The formula (\ref{ssmon}) is the result of
the computation this character in the Gelfand-Tsetlin basis.

\section{Semiinfinite forms and vertex operators}
\subsection{Semiinfinite forms}
Let the space $V$ be spanned by $\ul{k}$, $k\in\Z+\oh$.
The space $F=\LV$  is, by definition, spanned by vectors
$$
v_S=\ul{s_1} \wedge \ul{s_2} \wedge  \ul{s_3} \wedge  \dots\,,
$$
where $S=\{s_1>s_2>\dots\}\subset \Z+\oh$ is such a subset that
both sets
$$
S_+ = S \setminus \left(\Z_{\le 0} - \oh\right) \,, \quad
S_- = \left(\Z_{\le 0} - \oh\right) \setminus S
$$
are finite. We equip $\LV$ with the inner product $(.  , . )$ in
which the basis $\{v_S\}$ is orthonormal. The space $F$ is also
called the fermionic Fock space.

The infinite Clifford algebra $Cl(V)$ is generated by elements
$\psi_k, \ \psi_k^*, \ k\in \Z+\oh$ with defining relations
\[
\psi_k\psi_l+\psi_l\psi_k=0, \ \psi_k^*\psi_l^*+\psi_l^*\psi_k^*=0, \
\psi_k\psi_l^*+\psi_l^*\psi_k=\delta_{k,l},
\]
It acts on the Fock space $F$ as

\[
\psi_k \left(\ul{s_1} \wedge \ul{s_2} \wedge  \ul{s_3} \wedge  \dots\right)
= \ul{k} \wedge \ul{s_1} \wedge \ul{s_2} \wedge  \ul{s_3} \wedge
\dots ,
\]
\begin{multline*}
\psi_k^*\left(\ul{s_1} \wedge \ul{s_2} \wedge \dots \ul{s_l}\wedge\ul{k}
\wedge\ul{s_{l+1}}\wedge  \dots \right) = \\
(-1)^l\ul{s_1} \wedge \ul{s_2}
\wedge \dots \ul{s_l}\wedge\ul{s_{l+1}}\wedge  \dots  ,
\end{multline*}
\[
\psi_k^* \left(\ul{s_1} \wedge \ul{s_2} \wedge  \ul{s_3} \wedge
\dots\right)= 0, k\in \Z\backslash S
\]

The space $F$ is an irreducible representation of $Cl(V)$.
Notice that the operator representing $\psi_k^*$ is conjugate to the operator
representing $\psi_k$ with the respect to the scalar product in $F$.

The Lie algebra $\gl_\infty$ of $\Z\times \Z$-matrices with
finitely many entries acts naturally on $V$ and therefore acts
diagonally on semi-infinite wedge space $F$. This action extends
to the action of $a_\infty$ \cite{Kac} and is reducible. Irreducible
components $F(m)$ are eigenspaces of the operator
\[
C=\sum_{k<0}\psi_k^*\psi_k-\sum_{k<0}\psi_k\psi_k^*
\]
The Fock space decomposes into the direct sum
\[
F=\oplus_{m\in \Z} F(m)
\]
of irreducible representations of $a_\infty$ \cite{Kac}.

The subspace $F(m)$ is spanned by vectors $\ul{s_1} \wedge
\ul{s_2} \wedge \ul{s_3} \wedge  \dots$ with $s_i=m-i+1/2$ for
sufficiently large $i$. It is generated by the action of
$a_\infty$ on the vacuum vector ( the $a_\infty$ highest weight
vector in $F^{(m)}$ ):
\[
v^{(m)}_0=\ul{m-1/2}\wedge\ul{m-3/2}\wedge\ul{m-5/2}\wedge\dots
\]

The operators $\psi_k$ and $\psi_k^*$ act on $v^{(m)}_0$ as
\begin{eqnarray}
\psi_kv^{(m)}_0 &=&0, \ k\le m-1/2 , \\
\psi_k^*v^{(m)}_0&=&0, \ k > m-1/2 ,
\end{eqnarray}
They can be regarded as $a_\infty$-intertwining operators $\psi:
V\otimes F^{(m)}\to F^{(m+1)}$. where $V={\mathbb C}^{\mathbb
Z+\oh}$ is the vector representation of $\gl_\infty$.

Vectors in the space $F^{(m)}$ can be parameterized by partitions.
For a partition $\lambda$ define
\[
v^{(m)}_\lambda=\ul{\lambda_1+m-1/2}\wedge\ul{\lambda_2+m-3/2}\wedge\dots
\]
It is clear that these vectors span $F^{(m)}$.

\subsection{Clifford algebra and vertex operators}\label{vo}

Consider elements
$$
\al_n = \sum_{k\in\Z+\frac12} \psi_{k+n}\, \psi^*_k \,,
\quad n=\pm 1,\pm2, \dots \,,
$$

They satisfy commutation relations
\begin{eqnarray}\label{cmm}
\left[\al_n, \al_m\right] &=& -n \, \delta_{n,-m} , \\
\left[\alpha_n,\psi_k\right]&=&\psi_{k+n} , \\
\left[\alpha_n,\psi_k^*\right]&=-&\psi_{k-n}^*
\end{eqnarray}

It is clear that
\[
\alpha_n v^{(m)}_o=0
\]
for $n<0$ and $m\in \Z$.
Vertex operators are the formal power series
\[
\Gamma_+(x)=\exp(\sum_{n\ge 1}\frac{x^n}{n}\alpha_n), \
\Gamma_-(x)=\exp(\sum_{n\ge 1}\frac{x^n}{n}\alpha_{-n})
\]
The operator $\Gamma_-(x)$ acts finitely
in the space $F$ ( applied to any vector
of $F$ it acts as a polynomial in $x$). In particular
\[
\Gamma_-(x)v^{(m)}_0=v^{(m)}_0
\]

The operator $\Gamma_+(x)$ is conjugate to $\Gamma_-(x)$:
\begin{equation}\label{Gpm}
(\Gamma_-(x)v,w)=(v,\Gamma_+(x)w)
\end{equation}
and since the scalar product is symmetric
\[
(\Gamma_-(x)v,w)=(\Gamma_+(x)w,v)
\]

Notice that its action is defined in $F$
not only as a formal power series in $x$.
In a weak sense operators $\Gamma_\pm$ are operator-
valued functions which are analytic at $x=0$.

Define the formal Fourier transform of $\psi_k,\ \psi_k^*$ as
power series
\begin{equation}
\psi(z)=\sum_{k\in\Z+1/2} z^k \, \psi_k \,,
\quad
\psi^*(z)=\sum_{k\in\Z+1/2} z^{-k} \, \psi^*_k \,.\label{psiz}
\end{equation}

These operators and vertex operators satisfy the following
commutation relations:
\begin{equation} \label{gg}
\Gamma_+(x)\Gamma_-(y)=(1-xy)\Gamma_-(y)\Gamma_+(x),
\end{equation}
\begin{eqnarray}\label{crel}
\Gamma_+(x)\psi(z)&=&(1-z^{-1}x)^{-1}\psi(z)\Gamma_+(x) \\
\Gamma_-(x)\psi(z)&=&(1-xz)^{-1}\psi(z)\Gamma_-(x) \\
\Gamma_+(x)\psi^*(z)&=&(1-z^{-1}x)\psi^*(z)\Gamma_+(x) \\
\Gamma_-(x)\psi^*(z)&=&(1-xz)\psi^*(z)\Gamma_-(x)
\end{eqnarray}
Here left and right sides are corresponding formal power series.

The following is a well known statement.
\begin{theorem} The following identity holds:
\begin{equation}\label{me}
(\Gamma_-(x_1)\Gamma_-(x_2)\dots\Gamma_-(x_n)v^{(m)}_\lambda,
v^{(m)}_\mu)=
s_{\lambda/\mu}(x_1,\dots,x_n)
\end{equation}
\end{theorem}

This is a well known statement, for a proof see for example \cite{Kac}.
we will give a proof of it here. The key step is to show that the
identity (\ref{me}) holds for $n=1$ which can be easily derived from
the fact that
\[
v_\lambda^{(m)}=\psi_{j_1}\dots\psi_{j_n}v_0^{(m-n)}
\]
where $j_1=\lambda_1+m-1/2,j_2=\lambda_2+m-3/2,\dots$.
and from the identity
\[
\Gamma_-(x)\psi(z_1)\dots\psi(z_n)v_0^{(m-n)}=\prod_{i=1}^n(1-xz_i)^{-1}
\psi(z_1)\dots\psi(z_n)v_0^{(m-n)}
\]
for generating functions.

For matrix elements of generating functions $\psi(z)$ and $\psi^*(z)$
we have:
\begin{equation}\label{pmel}
\begin{array}{lc}
(\psi^*(w)\psi(z)v^{(m)}_0,v^{(m)}_0)=\left(\frac{z}{w}\right)^{m+1/2}
\frac{1}{1-z/w} & , \  |z|<|w| \\
(\psi(z)\psi^*(w)v^{(m)}_0,v^{(m)}_0)=\left(\frac{z}{w}\right)^{m-1/2}
\frac{1}{1-w/z} & , \  |z|>|w| .
\end{array}
\end{equation}

These identities follows from the summation of a geometric
series.

Let $A_i^a=\sum_jA^a_{ij}\psi_j$ and $B_i^a=\sum_j
B^a_{ij}\psi_j^*$ for $a=1,\dots, n$. The following identity is
known as a Wick's lemma:

\begin{equation}\label{wl}
(A^1_{i_1}B^1_{j_1}A^2_{i_2}B^2_{j_2}\dots A^n_{i_n}B^1_{j_n}
v^{(m)}_0,v^{(m)}_0)=\det( K_{ab})_{1\leq a,b\leq n}
\end{equation}
where
\[
K_{ab}=\left\{\begin{array}{ll}
(A^a_{i_a}B^b_{j_b},v^{(m)}_0,v^{(m)}_0), & \ \ a\leq b \\
-(B^b_{i_b}A^a_{j_a},v^{(m)}_0,v^{(m)}_0), & \ \ a> b
\end{array} \right.
\]

\section{Some asymptotic for limit shapes}

\subsection{}

The frozen boundary is singular at $\tau=V_j$.
When $\tau\to V_j$ the singular branch of the boundary curve
behave as
\begin{multline*}
\chi(\tau)=2\ln(\frac{|\tau-V_j|}{2})+|V_j|/2+\sum_{i=0}^{j-1}
(V_{i+1}-U_i)+ \\ \ln(\frac{\prod_{k\neq j, 1\leq k\leq
N}|1-e^{-U_k+U_j}|}{\prod_{1\leq k\leq N}|1-e^{-U_k+V_j}|})
+O(|\tau-V_j|) \ .
\end{multline*}

\subsection{}The limit curve is tangent to the lines $\tau=U_0$ and $\tau=U_N$
at the points $(\chi_L,U_0)$ and $(\chi_R,U_N)$ where
\[
\chi_L=\sum_{j=1}^N\ln(\frac{1-e^{U_0-V_j}}{1-e^{U_0-U_j}})+\frac{U_0}{2}
\]
\[
\chi_R=\sum_{j=1}^N\ln(\frac{1-e^{U_N-V_j}}{1-e^{U_N-U_j}})+U_0-\frac{U_N}{2}
\]

If $\tau=U_0+\epsilon$ or $\tau=U_N-\epsilon$ and
$\epsilon\to +0$ the asymptotic of two double critical points
$z^{(\pm)}(\tau)$ is
\[
z^{(\pm)}(\tau)=\left\{\begin{array}{ll}
e^{U_0}(1 \pm \sqrt{\epsilon}\sqrt C_L
+O(\epsilon)) \mbox{ if $\tau=U_0+\epsilon$}\\
e^{U_N}(1 \pm \sqrt{\epsilon}\sqrt C_R
+O(\epsilon)) \mbox{ if $\tau=U_N-\epsilon$}
\end{array}\right.
\]
where $C_L$ and $C_R$ are some constants.

The boundary curve near this point behave as:
\[
\chi(\tau)=\left\{\begin{array}{ll}
\chi_L(1+\sqrt{(\tau-U_0}D_L+O(|U_0-\tau)|,  \mbox{ if $\tau\to U_0+0$}\\
\chi_L(1+\sqrt{(U_N-\tau}D_R+O(|U_N-\tau|),  \mbox{ if $\tau\to U_N-0$}
\end{array}\right.
\]
where, again, $D_L$ and $D_R$ are some constants.

Similar solution exists near each point $\tau=U_j$, $j=1,\dots,
N-1$.

\subsection{}

Let us find the asymptotic of the density function near the
left boundary of the limit shape. Assume
\[
\chi=\chi_L+\delta\chi, \ \tau=U_0+\delta\tau
\]
for some positive $\delta\tau\to 0$.

Solutions to (\ref{cp2}) have the asymptotic $z=\exp(U_0)(1+\delta
z)$. Let us find $\delta z$ as a function of $\delta\chi$ and
$\delta\tau$. We have the following asymptotical expansions:
\[
\prod_{j=1}^N\frac{1-ze^{-V_j}}{1-ze^{-U_j}}=\prod_{j=1}^N\frac{1-e^{U_0-V_j}}{1-e^{U_0-U_j}}
(1-C\delta z +O({\delta z}^2))
\]
where
\[
C=\sum_{i=1}^N(\frac{e^{U_0-V_j}}{1-e^{U_0-U_j}}-\frac{e^{U_0-U_j}}{1-e^{U_0-U_j}}>0
\]
Keeping leading orders in $\delta\chi,\delta\tau$, and $\delta z$
in (\ref{cp2}) we obtain the equation for $\delta z$:
\[
C{\delta z}^2+\delta z((C+1/2)\delta\tau-\delta\chi)+\delta\tau=0
\]
In the region $\delta\tau\propto(\delta\chi)^2$ we have two
asymptotical solutions
\[
\delta z_{1,2}\simeq
\frac{\delta\chi}{2}\pm\sqrt{\frac{\delta\chi^2}{2}-\frac{\delta\tau}{C}}
\]
From here we obtain the asymptotic of the density function in this
region:
\[
\rho=\frac{\theta}{\pi}\simeq\arctan(\sqrt{\frac{4\delta\tau}{C\delta\chi^2}-1})
\]
Here $\theta$ is the argument of $\delta z_1$.
Notice that $\rho\to 1/2-0$ as $\delta\chi\to +0$. 

\section{The symmetry of correlation functions}

Change variables in the integral
\begin{multline}
K((t_1,h_1),(t_2,h_2))=\\ \frac{1}{(2\pi i)^2}
\int_{|z|<min\{1,R(t_1)}\int_{R^*(t_2)<|w|<1}
\frac{\Phi_-(z,t_1)\Phi_+(w,t_2)}{\Phi_+(z,t_1)\Phi_-(w,t_2)}
\frac{\sqrt{zw}}{z-w}z^{-j_1}w^{j_2}\frac{dzdw}{zw}
\end{multline}
from $z$ to $w^{-1}$ and from $w$ to $z^{-1}$. It becomes
\begin{multline}
K((t_1,h_1),(t_2,h_2))=\\ \frac{1}{(2\pi i)^2}
\int_{|w|^{-1}<min\{1,R(t_1)}\int_{R^*(t_2)<|z|^{-1}<1}
\frac{\tilde{\Phi}_-(z,-t_2)\tilde{\Phi}_+(w,t_1)}{\tilde{\Phi}_+(z,t_2)
\tilde{\Phi}_-(w,t_1)}
\frac{\sqrt{zw}}{z-w}z^{-j_2}w^{j_1}\frac{dzdw}{zw}
\end{multline}
where
\[
\tilde{\Phi}_+(z,t)=\prod_m(1-z\tilde{x}^+_m)
\]
\[
\tilde{\Phi}_-(z,t)=\prod_m(1-z^-\tilde{x}^-_m)
\]
and $\tilde{x}^\pm_m=x^\mp_{-m}$.

Thus, we have the following ``reflection'' symmetry of correlation
functions:
\[
K((j_1,t_1),(j_2,t_2))=\tilde{K}((j_2,-t_1)(j_1,-t_1)) \ ,
\]

This symmetry is obvious on the ``microscopical level''.
It corresponds to the reflection of the tiling in $t$-direction.

\end{document}